\documentclass{article}
\usepackage{latexsym,amssymb,lastpage}
\usepackage{graphicx,amsfonts}
\usepackage{times,mathptmx,bm,amsmath}
\usepackage{dcolumn}

\numberwithin{equation}{section}

\usepackage{a4wide}

\usepackage{amsfonts}
\usepackage{amsmath}
\usepackage{amssymb}
\usepackage{enumerate}
\usepackage{framed}
\usepackage{graphicx}
\usepackage{lscape}
\usepackage{bm}
\usepackage{color}
\usepackage{multicol}

\usepackage{xspace}
\usepackage{algorithm}
\usepackage{subfigure}
\usepackage{amsfonts}
\usepackage{amsmath}
\usepackage{enumerate}
\usepackage{framed}

\newtheorem{theorem}{Theorem}[section]

\newtheorem{lemma}[theorem]{Lemma}
\newtheorem{proposition}[theorem]{Proposition}

\newtheorem{assumption}[theorem]{Assumption}
\newtheorem{remark}[theorem]{Remark}

\newcommand{\figref}[1]{{Figure~\ref{#1}}}

\newcommand{\HH}{\mathbb{H}}

\newcommand{\N}{\mathbb{N}}

\usepackage{latexsym,amssymb,lastpage}
\usepackage{graphicx,amsfonts}
\usepackage{times,mathptmx,bm,amsmath}
\usepackage{dcolumn}
\usepackage{xspace}
\newcommand{\Leja}{{L\'{e}ja}\xspace}

\newcommand{\Dt}{\Delta t}

\newcommand{\thmref}[1]{{Theorem~\ref{#1}}}
\newcommand{\lemref}[1]{{Lemma~\ref{#1}}}
\newcommand{\secref}[1]{{Section~\ref{#1}}}
\newcommand{\assref}[1]{{Assumption~\ref{#1}}}
\newcommand{\propref}[1]{{Proposition~\ref{#1}}}

\begin{document}
\title{STOCHASTIC EXPONENTIAL INTEGRATORS FOR FINITE ELEMENT DISCRETIZATION OF SPDEs FOR MULTIPLICATIVE  \& ADDITIVE NOISE}
\author{{\sc Gabriel J Lord\footnote{G.J.Lord@hw.ac.uk} and Antoine Tambue \footnote{ Corresponding author: at150@hw.ac.uk (tambuea@gmail.com)}}\\[2pt]
Department of Mathematics and Maxwell Institute, \\[6pt]
 Heriot-Watt University, Edinburgh EH14 4AS, UK.\\[6pt]}
   

\pagestyle{headings}
\markboth{G.J.Lord and A.Tambue}{\rm Stochastic Exponential Integrators for finite
  element discretization of  SPDEs with general trace class noise. }
\maketitle
\begin{abstract}
{We consider the numerical approximation of a general second order
 semi--linear parabolic stochastic partial differential equation
 (SPDEs) driven by space-time noise, for multiplicative and additive noise. 
We examine convergence of exponential integrators for multiplicative and additive noise.
 We consider  noise that is in trace class and give a convergence proof in the root mean square $L^{2}$ norm.
We discretize in space with the finite element method and in our implementation we examine both the finite element 
and the finite volume methods.
We present results for a linear  reaction diffusion equation  
 in two dimensions as well as a nonlinear example of two-dimensional
 stochastic advection diffusion reaction equation motivated from
 realistic porous media flow.} 
{Parabolic stochastic partial differential equation, finite element, 
  exponential integrators, strong numerical approximation, multiplicative noise, additive noise.}
\end{abstract}
%

\section{Introduction}
We analyse the strong numerical approximation of an Ito
stochastic partial differential equation  defined in $\Omega\subset
\mathbb{R}^{d}$. Boundary conditions on the domain $\Omega$
are typically Neumann, Dirichlet or some mixed conditions.
We consider 
\begin{eqnarray}
  \label{adr}
  dX=(AX +F(X))dt + B(X)d W, \qquad  X(0)=X_{0},\qquad t \in [0, T],\qquad T>0
\end{eqnarray}
in a Hilbert space $H$. Here 
$A$ is the generator of an analytic semigroup $S(t):=e^{t A}, t\geq 0$
not neccessary self adjoint. The functions $F$ and $B$  are nonlinear of 
$X$ and the noise term $W(x,t)$ is a $Q$-Wiener process  defined on a filtered probability space $(\mathbb{D},\mathcal{F},\mathbf{P},\left\lbrace F_{t}\right\rbrace_{t\geq 0})$, that is 
white in time. 
The noise can be represented as a series in the eigenfunctions of the covariance
operator $Q$  given by
\begin{eqnarray}
  \label{eq:W}
  W(x,t)=\underset{i \in
    \mathbb{N}^{d}}{\sum}\sqrt{q_{i}}e_{i}(x)\beta_{i}(t), 
\end{eqnarray}
where $(q_i,e_{i})$, $i\in \mathbb{N}^{d}$ are the eigenvalues and  eigenfunctions
of the covariance operator $Q$ and $\beta_{i}$ are
independent and identically distributed standard Brownian motions.
Precise assumptions on $A$, $F$,$B$ and $W$ are given in
\secref{scheme} and under these type of technical assumptions 
it is well known, see \cite{DaPZ,PrvtRcknr,Chw}  that the unique mild
solution of \eqref{adr} is given by 
\begin{eqnarray}
\label{eq1}
X(t)=S(t)X_{0}+\int_{0}^{t}S(t-s)F(X(s))ds +\int_{0}^{t}S(t-s)B(X(s))dW(s).
\end{eqnarray}
Typical examples of the above type of equation are stochastic
(advection) reaction diffusion equations 
arising, for example, in pattern formation in physics and mathematical biology.
We illustrate our work with  both a simple reaction diffusion
equation where we can construct an exact solution 
\begin{eqnarray}
 dX=\left(\nabla \cdot  \mathbf{D}\nabla X - \lambda X \right)dt+ dW
\end{eqnarray}
 as well as the stochastic advection reaction diffusion equation 
\begin{eqnarray}
 dX=\left(\nabla \cdot  \mathbf{D}\nabla X -  \nabla \cdot (\mathbf{q}
  X)-\frac{X}{\vert X \vert +1}\right)dt+ XdW,
\end{eqnarray}
where $ \mathbf{D}$ is the diffusion  tensor, $\mathbf{q}$ is the
Darcy velocity field \cite{sebastianb} and $\lambda$ is a constant depending of the reaction function. 
The study of numerical solutions of SPDEs is an active 
research area and there is an extensive literature on numerical
methods for SPDEs (\ref{adr}). For temporal discretizations the 
linear implicit Euler scheme is often used  \cite{shardlow05,Haus1}, 
spatial discretizations are usually achieved with finite element
 \cite{allen98:_finit, Stig1, Yn:04}, finite difference \cite{shardlow05,Haus1}  or spectral Galerkin method \cite{LR,Jentzen1,Jentzen2,Jentzen3,Jentzen4}.  

In the special case with additive noise, new schemes using linear 
functionals of the noise have recently been considered
\cite{Jentzen1,Jentzen2,Jentzen3,Jentzen4,GTambue,GTambueexpo}.
 The finite element method is used for the spatial discretization in \cite{GTambue,GTambueexpo} and the spectral
Galerkin in  \cite{Jentzen1,Jentzen2,Jentzen3,Jentzen4}. 
Our schemes here are based on  using the finite element method (or finite volume method) for space discretization
so that we gain the  flexibility of
theses methods to deal with complex boundary
conditions and we can apply well developed techniques such as
upwinding to deal with advection. One of our schemes is the non--diagonal version of the stochastic scheme presented in \cite{LR,KLNS}
and the other is the extension of the deterministic exponential time differencing of order one \cite{MC} to stochastic exponential scheme. 
 Comparing to the schemes presented in \cite{GTambue,GTambueexpo} on additive noise, the results here are more general since
the linear operator $A$ does not need to be
 self adjoint and we do not need information about eigenvalues and eigenfunctions of the linear operator $A$. 
Furthermore we examine here convergence for Ito multiplicative noise for the 
exponential integrators, which has not so far been considered for SPDEs for these integrators.
As in \cite{GTambueexpo},
schemes presented here are based on exponential  matrix computation,
which is a notorious problem in numerical analysis
\cite{CMCVL}. However, new developments for both \Leja points and
Krylov subspace techniques \cite{kry,SID,Antoine,LE2,LE1,LE} have led
to efficient methods for computing matrix  exponentials.  
The convergence proof given below is similar to one in \cite{Yn:05}
for a finite element discretization in space and backward Euler based
method in time.  
The paper is organised as follows. In \secref{scheme} we
present the two  numerical schemes based on the exponential
integrators and our assumptions on \eqref{adr}. We also present and comment
on our convergence results. \secref{proof} contains
the proofs of our convergence theorems. 
We conclude in \secref{simulation} by presenting some simulations and
discuss implementation of these methods. 

\section{Mild solution, numerical schemes and main results}
\label{scheme}
\subsection{The abstract setting and mild solution}
Let us start by presenting briefly the notation for the main function
spaces and norms that we use in the paper. 
We denote by $\Vert \cdot \Vert$ the norm associated to
the inner product $(\cdot ,\cdot )$ of the Hilbert space $H$.
For a Banach space $\mathcal{V}$ we denote by $\Vert.\Vert_{\mathcal{V}}$ the norm of $\mathcal{V}$,
$L(\mathcal{V})$  the set of bounded linear mapping  from
$\mathcal{V}$ to $\mathcal{V}$
and  by $L_{2}(\mathbb{D},\mathcal{V})$ the space defined by
\begin{eqnarray}
\label{l2}
 L_{2}(\mathbb{D},\mathcal{V}):=\left\lbrace v\; \text{random variable with value in}\;\mathcal{V}:\;\mathbf{E} \Vert v\Vert_{\mathcal{V}}^{2}= 
\int_{\mathbb{D}} \Vert v(\omega)\Vert_{\mathcal{V}}^{2}d\mathbf{P}(\omega)<\infty \right\rbrace.
\end{eqnarray}
Let $Q: H\rightarrow H$ be a trace class operator.
 We introduce the spaces and
notation we need to define the $Q$-Wiener process. 
An operator $T \in L(H)$ is Hilbert-Schmidt if
\begin{eqnarray*}
 \Vert T\Vert_{HS}^{2}:=\underset{i\in \mathbb{N}^{d}}{\sum}\Vert T e_{i}
 \Vert ^2 < \infty,
\end{eqnarray*}
where $ (e_{i})$ is an orthonormal basis in H.
The sum  in $\Vert .\Vert_{HS}^{2}$ is independent of the choice of
the orthonormal basis in $H$. 
We denote the space of Hilbert--Schmidt operators from 
$Q^{1/2}(H)$ to $H$   by $L_{2}^{0}:=  HS(Q^{1/2}(H),H)$ 
and the corresponding norm $\Vert . \Vert_{L_{2}^{0}}$ by
\begin{eqnarray*}
 \Vert T\,\Vert_{L_{2}^{0}} := \Vert T
 Q^{1/2}\Vert_{HS}=\left( \underset{i \in \mathbb{N}^{d}}{\sum}\Vert
   T Q^{1/2} e_{i} \Vert^{2}\right)^{1/2},\qquad \ T\in L_{2}^{0}.  
\end{eqnarray*}
Let  $\varphi $  be a $L_{2}^{0}-$ process,  we have the following
equality using the Ito's isometry  \cite{DaPZ} 
\begin{eqnarray*}
 \mathbf{E} \Vert \int_{0}^{t}\varphi dW \Vert^{2}=\int_{0}^{t}
 \mathbf{E} \Vert \varphi \Vert_{L_{2}^{0}}^{2}ds=\int_{0}^{t}
 \mathbf{E} \Vert \varphi Q^{1/2} \Vert^{2}_{HS}ds,\qquad \qquad t\in [0,T]. 
\end{eqnarray*}
Let us give some assumptions required both for the existence and
uniqueness of the solution of equation (\ref{adr}) and for our
convergence proofs below.
\begin{assumption}
\label{assumptionn}
The operator $A$ is the generator of an analytic semigroup $S(t)=e^{t A}, \quad t\geq 0$.
\end{assumption}
In the Banach space  $\mathcal{D}((-A)^{\alpha/2})$, $\alpha \in
\mathbb{R}$, we use the notation 
$ \Vert (-A)^{\alpha/2}. \Vert =:\Vert .\Vert_{\alpha} $.
We recall some basic properties of the semigroup $S(t)$ generated by $A$.
\begin{proposition}
 \textbf{[Smoothing properties of the semigroup \cite{Henry}]}\\
\label{prop1}
Let $ \alpha >0,\;\beta \geq 0 $ and $0 \leq \gamma \leq 1$, then  there exist  $C>0$ such that
\begin{eqnarray*}
 \Vert (-A)^{\beta}S(t)\Vert_{L(H)} &\leq& C t^{-\beta}\;\;\;\;\; \text {for }\;\;\; t>0\\
  \Vert (-A)^{-\gamma}( \text{I}-S(t))\Vert_{L(H)} &\leq& C t^{\gamma} \;\;\;\;\; \text {for }\;\;\; t\geq 0.
\end{eqnarray*}
In addition,
\begin{eqnarray*}
(-A)^{\beta}S(t)&=& S(t)(-A)^{\beta}\quad \text{on}\quad \mathcal{D}((-A)^{\beta} )\\
\text{If}\;\;\; \beta &\geq& \gamma \quad \text{then}\quad
\mathcal{D}((-A)^{\beta} )\subset \mathcal{D}((-A)^{\gamma} ),\\
\Vert D_{t}^{l}S(t)v\Vert_{\beta}&\leq& C t^{-l-(\beta-\alpha)/2} \,\Vert v\Vert_{\alpha},\;\; t>0,\;v\in  \mathcal{D}((-A)^{\alpha/2})\;\; l=0,1, 
\end{eqnarray*}
where $ D_{t}^{l}:=\dfrac{d^{l}}{d t^{l}}$.
\end{proposition}
We describe now in detail the assumptions that we make
on the nonlinear terms $F$,$B$ and the noise $W$.
\begin{assumption}
\label{assumption1}
\textbf{[Assumption on the drift term $F$]}
There exists  a positive constant $L> 0$  such  that
$F$ is continuous in $H$  and satisfies the following Lipschitz condition
\begin{eqnarray*}
 \Vert F(Z)- F(Y)\Vert \leq L \Vert Z- Y\Vert \qquad \forall \quad Z, \;  Y  \in H. 
\end{eqnarray*}
\end{assumption}
As a consequence, there  exists a constant $C>0$  such that
\begin{eqnarray*}
 \Vert F(Z) \Vert &\leq&   \Vert F (0)\Vert + \Vert F (Z) - F (0)\Vert
 \leq \Vert F (0)\Vert + L \Vert Z\Vert \leq C( 1
 +\Vert Z \Vert )\qquad \qquad Z\in H.\\
\end{eqnarray*}
\begin{assumption}
 \label{assumption2}
\textbf{[Assumption on the noise and the diffusion term $B$]}\\
The covariance operator $Q$ is in reace class i.e. Tr($Q)< \infty $, and  there exists  a positive constant $L> 0$ such  that
 $B$ is continuous in $H$  and satisfies the following condition
\begin{eqnarray*}
 \Vert B(Z)- B(Y)\Vert_{L_{2}^{0}} \leq L \Vert Z- Y\Vert \qquad \forall Z,
 Y  \in H.
\end{eqnarray*}
\end{assumption}
As a consequence, there exists a constant  $C>0$  such that
\begin{eqnarray*}
 \Vert B(Z) \Vert_{L_{2}^{0}} &\leq&   \Vert B(0)\Vert_{L_{2}^{0}} + \Vert B(Z) - B(0)\Vert_{L_{2}^{0}}
 \leq \Vert B(0)\Vert_{L_{2}^{0}} + L \Vert Z\Vert \leq C( 1
 +\Vert Z \Vert )\qquad \qquad Z\in H.\\
\end{eqnarray*}

\begin{theorem}
 \label{existth}
 \textbf{[Existence and uniqueness \cite{DaPZ}]}\\
Assume that the initial solution $X_{0}$ is an $F_{0}-$measurable $H-$valued random variable and \assref{assumption1}, 
\assref{assumption2} are satisfied.  
There exists a mild solution $X$ to \eqref{adr} unique, up to
equivalence among the processes, satisfying 
\begin{eqnarray}
 P\left( \int_{0}^{T}\Vert X(s)\Vert^{2}ds < \infty\right).
\end{eqnarray}
 For any $p\geq 2$ there exists a constant $C =C(p,T)>0 $ such that 
\begin{eqnarray}
\label{ineq2}
 \underset{t\in [0,T]}{\sup}\mathbf{E} \Vert X(t)\Vert^{p} \leq C\left(1+\mathbf{E} \Vert X_{0} \Vert^{p}\right). 
\end{eqnarray}
 For any $p>2$ there exists a constant $C_{1} =C_{1}(p,T)>0 $ such that 
\begin{eqnarray}
 \mathbf{E}\underset{t\in [0,T]}{\sup} \Vert X(t)\Vert^{p} \leq C_{1}\left(1+\mathbf{E} \Vert X_{0} \Vert^{p}\right). 
\end{eqnarray}
\end{theorem}
The following theorem proves a regularity result of the mild solution $X$ of \eqref{adr}.
\begin{theorem}
 \label{newtheo}
Assume that  \assref{assumption1} and \assref{assumption2} hold. Let $X$ be the  mild solution of (\ref{adr}) given in (\ref{eq1}). If  $X_{0} \in L_{2}(\mathbb{D},\mathcal{D}((-A)^{\beta/2})),\, \beta \in [0,1)$ 
then  for all  $ t\in [0,T],\,X(t) \in L_{2}(\mathbb{D},\mathcal{D}((-A)^{\beta/2}))$  with 
\begin{eqnarray*}
\label{regsoluion}
 \left(\mathbf{E}\Vert X(t) \Vert_{\beta}^{2}\right)^{1/2}\leq C \left(1+\left(\mathbf{E}\Vert X_{0}\Vert_{\beta}^{2}\right)^{1/2}+\left(1+\mathbf{E} \Vert X_{0} \Vert^{2}\right)^{1/2}\right).
\end{eqnarray*}
\end{theorem}

{\bf Proof}
Recall  that if $X$ is the mild solution of (\ref{adr}), according to \eqref{l2} we need to estimate $ \left(\mathbf{E}\Vert X(t) \Vert_{\beta}^{2}\right)^{1/2}$ and check that
$$
\left(\mathbf{E}\Vert X(t) \Vert_{\beta}^{2}\right)^{1/2}<\infty.
$$

Recall that the mild solution is given by 
$$X(t)=S(t)X_{0}+\int_{0}^{t}S(t-s)F(X(s))ds
+\int_{0}^{t}S(t-s)B(X(s))dW(s)$$ 
then
\begin{eqnarray*}
 \left(\mathbf{E}\Vert X(t) \Vert_{\beta}^{2}\right)^{1/2} &\leq & \left(\mathbf{E}\Vert S(t)X_{0}\Vert_{\beta}^{2}\right)^{1/2} +\left(\mathbf{E}\Vert\int_{0}^{t}S(t-s)F(X(s))ds\Vert_{\beta}^{2}\right)^{1/2}+ \left(\mathbf{E}\Vert\int_{0}^{t}S(t-s)B(X(s))dW(s) \Vert_{\beta}^{2}\right)^{1/2}\\
                              &=& I+II+III.  
\end{eqnarray*}
Since $X_{0} \in L_{2}(\mathbb{D},\mathcal{D}((-A)^{\beta/2})),\, \beta \in [0,1)$, we obviouly have 
$$
I=\left(\mathbf{E}\Vert S(t)X_{0}\Vert_{\beta}^{2}\right)^{1/2} \leq C \,\left(\mathbf{E}\Vert X_{0}\Vert_{\beta}^{2}\right)^{1/2}.
$$
As a  consequence of \assref{assumption1} and the semigroup properties in \propref{prop1} we have 
\begin{eqnarray*}
II=\left(\mathbf{E}\Vert\int_{0}^{t}S(t-s)F(X(s))ds\Vert_{\beta}^{2}\right)^{1/2}&\leq& \int_{0}^{t}\left(\mathbf{E}\Vert S(t-s)F(X(s))\Vert_{\beta}^{2}\right)^{1/2}ds\\
  &\leq& \int_{0}^{t}\left(\mathbf{E}\Vert (-A)^{\beta/2}S(t-s) F(X(s))\Vert^{2} \right)^{1/2} ds\\
 &\leq& C \left(\int_{0}^{t}\Vert (-A)^{\beta/2}S(t-s)\Vert_{L(L^{2}(\Omega))} ds\right) \left(\underset{ 0\leq s\leq t}{\sup} \left(\mathbf{E}\left(1+\Vert X(s)\Vert\right)^{2} \right)^{1/2}\right)\\
&\leq& C \left(\int_{0}^{t}(t-s)^{-\frac{\beta}{2}}ds\right)\left(1+\underset{ 0\leq s\leq t}{\sup}\left(\mathbf{E}\Vert X(s)\Vert^{2}\right)^{1/2}\right)\\
 &\leq& C\,\left(1+\underset{ 0\leq s\leq t}{\sup}\left(\mathbf{E}\Vert X(s)\Vert^{2}\right)^{1/2}\right).
\end{eqnarray*}
Finally, Ito's isometry and the consequence of \assref{assumption2} yields
\begin{eqnarray*}
 III^{2}&=&\mathbf{E}\Vert\int_{0}^{t}S(t-s)B(X(s))dW(s) \Vert_{\beta}^{2}\\
   &=&\int_{0}^{t}\mathbf{E}\Vert (-A)^{\beta/2} S(t-s)B(X(s))\Vert_{L_{2}^{0}}^{2}ds\\
   &\leq & C \left(\int_{0}^{t}\Vert (-A)^{\beta/2}S(t-s)\Vert_{L(L^{2}(\Omega))}^{2} ds\right) \left(1+\underset{ 0\leq s\leq t}{\sup}\mathbf{E}\Vert X(s)\Vert^{2}\right)\\
&\leq & C \left(1+\underset{ 0\leq s\leq t}{\sup}\mathbf{E}\Vert X(s)\Vert^{2}\right),
\end{eqnarray*}
thus
\begin{eqnarray*}
 III \leq C \left(1+\left(\underset{ 0\leq s\leq t}{\sup}\mathbf{E}\Vert X(s)\Vert^{2}\right)^{1/2}\right).
\end{eqnarray*}
Then, using \eqref{ineq2} with $p=2$ yields
\begin{eqnarray*}
 \left(\mathbf{E}\Vert X(t) \Vert_{\beta}^{2}\right)^{1/2} &\leq& C \,\left(1+\left(\mathbf{E}\Vert X_{0}\Vert_{\beta}^{2}\right)^{1/2} +\underset{ 0\leq s\leq t}{\sup}\left(\mathbf{E}\Vert X(s)\Vert^{2}\right)^{1/2}\right)
\\
& \leq &  C \,\left(1+\left(\mathbf{E}\Vert X_{0}\Vert_{\beta}^{2}\right)^{1/2}+\left(1+\mathbf{E} \Vert X_{0} \Vert^{2}\right)^{1/2}\right),\\
& < & \infty.
\end{eqnarray*}
Then, if  $X_{0} \in L_{2}(\mathbb{D},\mathcal{D}((-A)^{\beta/2})),\, \beta \in [0,1)$ 
then  for all  $ t\in [0,T],\,X(t) \in L_{2}(\mathbb{D},\mathcal{D}((-A)^{\beta/2}))$.

More results about the regularity of the mild solution $X$ can be
found in \cite{MJentzen1,Mpato}. 

\subsection{Application to the second order  semi--linear parabolic SPDEs}
We assume that $\Omega$  has a smooth boundary or is a convex polygon of $\mathbb{R}^{d},\;d=1,2,3$. 
In the sequel of this paper, for convenience of presentation, we take $A$ to be a second order
operator as this simplifies the convergence proof.

More precisely we take $H=L^{2}(\Omega)$  and  consider the general second order 
semi--linear parabolic stochastic partial differential equation given by
\begin{eqnarray}
\label{sadr}
 dX(t,x)=\left(\nabla \cdot \textbf{D} \nabla X(t,x) -\mathbf{q} \cdot \nabla X(t,x) + f(x,X(t,x))\right) dt +b(x,X(t,x))dW(t,x),
\end{eqnarray}
$x \in \Omega, t\in[0,T]$ 
where $f,b:\Omega \times\mathbb{R} \rightarrow \mathbb{R}$ are
two continuously differentiable functions with globally bounded derivatives.

In the abstract form given in \eqref{adr}, the linear operator is defined by
\begin{eqnarray}
 A&=&\underset{i,j=1}{\sum^{d}}\dfrac{\partial }{\partial x_{i}}\left( D_{i,j}\dfrac{\partial 
 }{\partial x_{j}}\right) - \underset{i=1}{\sum^{d}}q_{i}\dfrac{\partial 
 }{\partial x_{i}},
\end{eqnarray}
where we assume that  $D_{i,j} \in L^{\infty}(\Omega), q_{i}\in L^{\infty}(\Omega)$ and that there exists a positive constant $c_{1}>0$ such that 
\begin{eqnarray}
\label{ellipticity}
\underset{i,j=1}{\sum^{d}}D_{i,j}(x)\xi_{i}\xi_{j}\geq c_{1}\vert \xi \vert^{2}  \;\;\;\;\;\;\forall \xi \in \mathbb{R}^{d}\;\;\; x \in \overline{\Omega}\;\;\; c_{1}>0,
\end{eqnarray}
and  $F : H \rightarrow H$, $ B : H\rightarrow HS(Q^{1/2}(H), H)$ are defined by
\begin{eqnarray}
\label{nemform}
(F(v))(x)=f(x,v(x)),\qquad (B(v)u)(x)=b(x,v(x))\cdot u(x),
\end{eqnarray}
for  all $ x\in \Omega,\;v\in H,\; u \in Q^{1/2}(H)$, with $H=L^{2}(\Omega)$. As the functions $f$ and $b$ are
two continuously differentiable functions with globally bounded derivatives, the  Nemytskii operator $F$ corresponding to $f$  and the multiplication operator $B$ defined in \eqref{nemform} satisfy  Assumptions \ref{assumption1}--\ref{assumption2} 
for appropriate eigenfunctions  such that $\underset{i \in \mathbb{N}^{d}}{\sup} \left[\underset{x\in \Omega}{\sup} \vert e_{i}(x)\vert\right]<\infty $ (see \cite[Section 4 ]{MJentzen1}).

Notice that  by the definitions of the operator $B$ and  $\Vert. \Vert_{L_{2}^{0}}$, for $Y  \in H=L^{2}(\Omega)$
\begin{eqnarray}
 \Vert B(Y)\Vert_{L_{2}^{0}}^{2}= \underset{ i\in \mathbb{N}^{d}}{\sum} \Vert b(Y)Q^{1/2}e_{i}\Vert^{2},
\end{eqnarray}
where $b(Y)$ is the Nemytskii operator defined by
\begin{eqnarray}
 b(Y)(x)=b(x,Y(x))\qquad\qquad x \in \Omega.
\end{eqnarray}

We introduce two spaces $\HH$ and $V$ where $\HH\subset V $ that depend
on the choice of the boundary conditions for the domain of the operator $A$ and the corresponding bilinear form.
For Dirichlet boundary conditions we let 
\begin{eqnarray*}
V= \HH= H_{0}^{1}(\Omega)=\{v\in H^{1}(\Omega): v=0\;\;
\text{on}\;\;\partial \Omega\}, 
\end{eqnarray*}
and for Robin boundary conditions, Neumann boundary being a special case, we take $V=  H^{1}(\Omega)$ and 
\begin{eqnarray*}
\HH = \left\lbrace v\in H^{1}(\Omega): \partial v/\partial
  \nu_{A}+\alpha_{0} v=0\quad \text{on}\quad \partial \Omega\right\rbrace, \qquad \alpha_{0} \in \mathbb{R}.
 \end{eqnarray*}
See \cite{lions} for details.
The corresponding bilinear form of $-A$ is given by
\begin{eqnarray}
\label{var}
a(u,v)=\int_{\Omega}\left(\underset{i,j=1}{\sum^{d}} D_{i,j}\dfrac{\partial u}{\partial x_{j}} \dfrac{\partial v}{\partial x_{i}}+\underset{i=1}{\sum^{d}}q_{i} \dfrac{\partial u}{\partial x_{j}}v\right)dx\;\;\;\;\;\;\; u, v \in V
\end{eqnarray}
for Dirichlet  and Neumann boundary conditions, and by
\begin{eqnarray}
\label{var1}
a(u,v)=\int_{\Omega}\left(\underset{i,j=1}{\sum^{d}} D_{i,j}\dfrac{\partial u}{\partial x_{j}} \dfrac{\partial v}{\partial x_{i}}+\underset{i=1}{\sum^{d}}q_{i} \dfrac{\partial u}{\partial x_{j}}v\right)dx
+\int_{\partial \Omega} \alpha_{0} u\,v\,dx \;\;\;\;\;\;\; u, v \in V,
\end{eqnarray}
for Robin boundary conditions.
According to G\aa{}rding's inequality (see \cite{ATthesis,lions}),
there exists two positive constants $c_{0}$ and $\lambda_{0}$  such
that 
 \begin{eqnarray}
 \label{coer}
  a(v,v)+c_{0}\Vert v\Vert^{2}\geq  \lambda_{0}\Vert v\Vert_{H^{1}(\Omega)}^{2}\;\;\; \quad \quad \forall v\in V.
\end{eqnarray}

By adding and subtracting $c_{0}X dt$ on the right hand side of (\ref{adr}), we have a new operator that we still call $A$  corresponding to the new bilinear form that we still call 
$a$ such that the following  coercivity property holds
\begin{eqnarray}
\label{ellip}
a(v,v)\geq \; \lambda_{0}\Vert v\Vert_{H^{1}(\Omega)}^{2}\;\;\;\;\;\forall v \in V.
\end{eqnarray}
Note that the expression of the nonlinear term $F$ has changed as we include the term $-c_{0}X$ 
in a new nonlinear term that we still denote by $F$.
The coercivity property (\ref{ellip}) implies that $A$ is a sectorial on $L^{2}(\Omega)$ i.e.  there exists $C_{1},\, \theta \in (\frac{1}{2}\pi,\pi)$ such that
\begin{eqnarray}
 \Vert (\lambda I -A )^{-1} \Vert_{L(L^{2}(\Omega))} \leq \dfrac{C_{1}}{\vert \lambda \vert }\;\quad \quad 
\lambda \in S_{\theta},
\end{eqnarray}
where $S_{\theta}=\left\lbrace  \lambda \in \mathbb{C} :  \lambda=\rho e^{i \phi},\; \rho>0,\;0\leq \vert \phi\vert \leq \theta \right\rbrace $ (see \cite{Henry,lions}).
 Then  $A$ is the infinitesimal generator of bounded analytic semigroups $S(t):=e^{t A}$  on $L^{2}(\Omega)$  such that
\begin{eqnarray}
S(t):= e^{t A}=\dfrac{1}{2 \pi i}\int_{\mathcal{C}} e^{ t\lambda}(\lambda I - A)^{-1}d \lambda,\;\;\;\;\;\;\;
 \;t>0,
\end{eqnarray}
where $\mathcal{C}$  denotes a path that surrounds the spectrum of $A $.

Functions in $\HH$ satisfy the boundary
conditions and with $\HH$ in hand we can characterize the
domain of the operator $(-A)^{r/2}$ and have the following norm
equivalence \cite{Stig,ElliottLarsson} for $r=1,2$
\begin{eqnarray*}
\Vert v \Vert_{H^{r}(\Omega)} \equiv \Vert (-A)^{r/2} v
\Vert=:\Vert  v \Vert_{r}, \qquad 
\forall v\in \mathcal{D}((-A)^{r/2})=  \HH\cap H^{r}(\Omega).
\end{eqnarray*} 
\subsection{Numerical schemes}
We consider the discretization of the spatial domain by a finite element 
triangulation.
Let $\mathcal{T}_{h}$ be a set of disjoint intervals  of $\Omega$ 
(for $d=1$), a triangulation of $\Omega$ (for $d=2$) or a set of
tetrahedra (for $d=3$) with maximal length $h$.
Let $V_{h}\subset V$ denote the space of continuous functions that are
piecewise linear over the triangulation $\mathcal{T}_{h}$.
To discretize in space we introduce the projection $P_h$ from $L^{2}(\Omega)$  onto $V_{h}$ defined  for  $u
\in L^{2}(\Omega)$ by 
\begin{eqnarray}
 (P_{h}u,\chi)=(u,\chi)\qquad \forall\;\chi \in V_{h}.
\end{eqnarray}
The discrete operator $A_{h}: V_{h}\rightarrow V_{h}$ is  defined by
\begin{eqnarray}
( A_{h}\varphi,\chi)=(A\varphi,\chi)=-a(\varphi,\chi)\qquad \varphi,\chi \in V_{h}.
\end{eqnarray}
Like the operator $A$, the discrete operator $A_h$ is also the generator of an analytic semigroup  $S_h:=e^{tA_{h}}$.

The semi--discrete in space version of the problem (\ref{adr}) is to
find the process $X^{h}(t)=X^{h}(.,t) \in V_{h}$ such  that for $t
\in[0, T]$,
\begin{eqnarray}
\label{dadr}
  dX^{h}=(A_{h}X^{h} +P_{h}F(X^{h}))dt + P_{h} B(X^{h}) d W,\qquad
  X^{h}(0)=P_{h}X_{0}.
\end{eqnarray}
The  mild solution of (\ref{dadr})  at  time $t_{m}=m \Delta t ,\;\;\Delta t>0\; $ is given by
\begin{eqnarray}
\label{dmild}
  X^{h}(t_{m})&=&S_{h}(t_{m})P_{h}X_{0}+\int_{0}^{t_{m}} S_{h}(t_{m}-s) P_{h}F(X^{h}(s))ds \nonumber \\ && +
 \int_{0}^{t_{m}} S_{h}(t_{m}-s)P_{h}B(X^{h})d W(s).
\end{eqnarray}
Then, given the mild solution at the time $t_{m}$,  we can construct the corresponding solution at $t_{m+1}$ as
\begin{eqnarray*}
 X^{h}(t_{m+1})&=&S_{h}(\Delta t)X^{h}(t_{m})+\int_{0}^{ \Delta t} S_{h}( \Delta t-s) P_{h}F(X^{h}(s+t_{m}))ds  \\ && +\int_{t_{m}}^{t_{m+1}} S_{h}(t_{m+1}-s)P_{h}B(X^{h})d W(s).
\end{eqnarray*}
To build the first  numerical scheme, we use the following approximations
\begin{eqnarray*}
 S_{h}( \Delta t-s)F(X^{h}( t_{m}+s))&\approx & S_{h}( \Delta t) F(X^{h}( t_{m}))\;\;\;\;\; \;\;\;\;\;s \in [0,\; \Delta t],\\
  S_{h}(t_{m+1}-s)P_{h}B(X^{h}(s))&\approx& S_{h}(\Delta t) P_{h}B(X^{h}(t_{m}))\;\;\; s\in [t_{m},t_{m+1}].
\end{eqnarray*}
We can define our approximation $Y_{m}^{h}$ of  $X(m \Delta t)$ by
\begin{eqnarray}
\label{new0}
 Y_{m+1}^{h}&=&e^{\Delta t A_{h}}\left(Y_{m}^{h}+P_{h}F(Y_{m}^{h})+ P_{h}B(Y_{m}^{h})\Delta W_{m}\right) \nonumber\\
 &=&\varphi_{0}(\Delta t A_{h})\left(Y_{m}^{h}+P_{h}F(Y_{m}^{h})+ P_{h}B(Y_{m}^{h})\Delta W_{m}\right),
\end{eqnarray}
where 
\begin{eqnarray*}
\varphi_{0}(\Delta t A_{h})&:=&e^{\Delta t A_{h}},\\
 \Delta W_{m} &:=& W_{m+1}-W_{m}=\sqrt{\Delta t} \underset{i \in \mathbb{N}^{d}}{\sum}\sqrt{q_{i}} R_{i,m}e_{i},
\end{eqnarray*}
with $ R_{i,m}$ are independent, standard normally distributed random
variables with means $0$ and variance $1$.
We call the scheme in (\ref{new0}) SETDM0. 
To build the second numerical scheme, we use the following approximations
\begin{eqnarray*}
 F(X^{h}( t_{m}+s))&\approx&  F(X^{h}( t_{m}))\;\;\;\;\; \;\;\;\;\;s \in [0,\; \Delta t],\\
 S_{h}(t_{m+1}-s)P_{h}B(X^{h}(s))&\approx& S_{h}(\Delta t) P_{h}B(X^{h}(t_{m}))\;\;\; s\in [t_{m},t_{m+1}]. 
\end{eqnarray*}
We can define our approximation $X_{m}^{h}$ of  $X(m \Delta t)$ by
\begin{eqnarray}
\label{new}
 X_{m+1}^{h}=e^{\Delta t A_{h}}X_{m}^{h}+A_{h}^{-1}\left( e^{\Delta t A_{h}}-I\right)P_{h}F(X_{m}^{h})+e^{\Delta t A_{h}}P_{h}B(X_{m}^{h})\left( W_{m+1}-W_{m}\right).
\end{eqnarray}
For efficiency we can rewrite the scheme (\ref{new}) as
\begin{eqnarray*}
 X_{m+1}^{h}=X_{m}^{h}+\Delta t \varphi_{1}(\Delta t A_{h})\left(A_{h} \left(X_{m}^{h}+P_{h}B(X_{m}^{h})\Delta W_{m}\right) +P_{h}F(X_{m}^{h})\right),
\end{eqnarray*}
where   $$\varphi_{1}(\Delta t A_{h})=(\Delta t\, A_{h})^{-1}\left( e^{\Delta t A_{h}}-I\right)= \frac{1}{\Delta t}\int_{0}^{\Delta t} e^{(\Delta t- s)A_{h}}ds. $$
We call the scheme in (\ref{new}) SETDM1. This scheme is also used in \cite{iyaboThesis} with the Fourier  method to solve fourth order stochastic problems.

\subsection{Main result}
Throughout the paper we
take $t_m=m\Dt \in (0,T]$, where $T=M\Dt$ for $m,M\in\N$. We take $C$
to be a constant that may depend on $T$ and other parameters but not
on $\Dt$ or $h$. 
Our  result is a strong convergence result in $L^2$ for schemes SETDM1 and SETDM0.
\begin{theorem}
\label{th1}
Let  $X(t_m)$  be the mild  solution of equation \eqref{adr}  at time \;$t_m=m\Delta t,\; \Delta t>0$ represented
by (\ref{eq1}). Let $\zeta_m^{h}$ be the numerical approximations through \eqref{new}  or \eqref{new0} 
($\zeta_m^{h}=X_m^{h} $ for scheme  SETDM1 and  $\zeta_m^{h}=Y_m^{h}$ for scheme
 SETDM0)  and  $ 0<\gamma< 1$. Assume that
$b(L_{2}(\mathbb{D},\mathcal{D}((-A)^{\alpha})))\subset L_{2}(\mathbb{D},\mathcal{D}((-A)^{\alpha}))$ for $\alpha \in (0,\gamma/10)$ small 
enough. The following estimates hold.

If  $X_{0} \in L_{2}(\mathbb{D},\mathcal{D}((-A)^{\gamma})), \; 0< \gamma \leq 1/2$ then
\begin{eqnarray*}
 \left(\mathbf{E}\Vert X(t_m)-\zeta_m^{h}\Vert^{2}\right)^{1/2}& \leq
 &C \left( t_{m}^{(-1+2 \gamma)/2}  
   h + \Delta t^{\gamma/2}\right).
\end{eqnarray*}
 If  $X_{0}\in L_{2}(\mathbb{D},\mathcal{D}((-A)^{\gamma})), \; 1/2 \leq \gamma < 1 $ then 
\begin{eqnarray*}
 \left(\mathbf{E}\Vert X(t_m)-\zeta_m^{h}\Vert^{2}\right)^{1/2}& \leq
 &C \left( h + \Delta t^{\gamma/2}\right).
\end{eqnarray*}
Suppose that  $  F(L_{2}(\mathbb{D},\mathcal{D}((-A)^{\alpha'})))\subset L_{2}(\mathbb{D},\mathcal{D}((-A)^{\alpha'})),\; \alpha' \in (0,\gamma/10)$ small enough:

If  $X_{0} \in L_{2}(\mathbb{D},\mathcal{D}((-A)^{\gamma}))$ then 
 \begin{eqnarray*}
 \left(\mathbf{E}\Vert X(t_m)-\zeta_m^{h} \Vert^{2}\right)^{1/2}& \leq
 &C \left(  t_{m}^{(-1+\gamma)} h^{2} + \Delta t^{\gamma/2}\right).
\end{eqnarray*}

If  $X_{0} \in L_{2}(\mathbb{D},\mathcal{D}(-A)) $ and  $ b(L_{2}(\mathbb{D},\mathcal{D}(-A)))\subset L_{2}(\mathbb{D},\mathcal{D}(-A)) $ then
\begin{eqnarray*}
 \left(\mathbf{E}\Vert X(t_m)-\zeta_m^{h}\Vert^{2}\right)^{1/2}& \leq
 &C \left(h^{2} + \Delta t^{(1/2-\epsilon)}\right).
\end{eqnarray*}
$\epsilon \in (0,1/2)$, small enough. 
\end{theorem}
\begin{remark}
 In the proof of  \thmref{th1}, the assumptions $b(L_{2}(\mathbb{D},\mathcal{D}((-A)^{\alpha})))\subset L_{2}(\mathbb{D},\mathcal{D}((-A)^{\alpha}))$ for $\alpha \in (0,\gamma/10)$ and
 $  F(L_{2}(\mathbb{D},\mathcal{D}((-A)^{\alpha'})))\subset L_{2}(\mathbb{D},\mathcal{D}((-A)^{\alpha'})),\; \alpha' \in (0,\gamma/10)$ are enough. We only need that
$b(X(t))\in  L_{2}(\mathbb{D},\mathcal{D}((-A)^{\alpha})),\; \alpha \in (0,\gamma/10)$ and  $F(X(t))\in  L_{2}(\mathbb{D},\mathcal{D}((-A)^{\alpha'})),\;\ \alpha'\in (0,\gamma/10),\;\; \forall \; t\in [0,T] $, where $X$ is the mild 
solution of equation \eqref{adr}.
\end{remark}

Although we have taken the linear operator $A$ to be a second order operator, similar results will hold, for higher 
order operators. Computationaly, the noise given by (\ref{eq:W}) is truncated to $N$ terms. Therefore the corresponding approximated solutions become $X_{m}^{h,N}$ for SETDM1 and 
$Y_{m}^{h,N}$ for scheme SETDM0. For noise where the eigenvalues of the covariance operator have a strong exponential decay, $X_{m}^{h,N}$ and $Y_{m}^{h,N}$ are close to $X_{m}^{h}$ and
 $Y_{m}^{h}$ respectively. 
In the case of additive noise, it has been proved in \cite{Stig1} that with the truncation to $N$ terms of the noise (\ref{eq:W}) the corresponding discrete mild solution  $X^{h,N}$ in  (\ref{dmild}) 
has the same order of accuracy respect to $h$ as $X^{h}$.

We also note that smooth noise improves the accuracy in  \thmref{th1} for additive noise (see \cite{KLNS} and \figref{FIG0022} in \secref{simulation}).

\section{Proofs of the main results} 
\label{proof}
\subsection{Some preparatory results}
We introduce  the Riesz representation operator  $R_{h}: V \rightarrow V_{h}$ defined by 
\begin{eqnarray*}
 (-A R_{h}v,\chi)=(-A v,\chi)=a(v,\chi)\qquad \qquad v \in V,\; \forall \chi \in V_{h}.
\end{eqnarray*}
Under the regularity assumptions on the triangulation and in view of $V-$ellipticity \eqref{ellipticity},  it is well known (see \cite{lions,ciarlet}) that the following error bounds holds
\begin{eqnarray}
\label{regulard}
  \Vert R_{h}v-v\Vert +h \Vert
    R_{h}v-v\Vert_{H^{1}(\Omega} \leq  C h^{r} \Vert
  v\Vert_{H^{r}(\Omega)}, \qquad v\in V\cap
  H^{r}(\Omega),\; \; r \in \{1,2\}. 
\end{eqnarray}

We start by examining the deterministic linear problem.
Find $u \in V$ such that such that  
\begin{eqnarray}
\label{homog}
u'=Au \qquad  
\text{given} \quad u(0)=v,\qquad  t\in (0,T] .
\end{eqnarray}
The corresponding semi-discretization in space is : Find $u_{h} \in
V_{h}$ such that  
$$u_{h}'=A_{h}u_{h}$$ 
where $u_{h}^{0}=P_{h}v$.
Define the operator 
\begin{eqnarray}
\label{form1}
T_{h}(t) :=  S(t)-S_{h}(t) P_{h} = e^{tA} - e^{tA_h}P_h
\end{eqnarray} 
so that $u(t)-u_{h}(t)= T_{h}(t) v$.

\begin{lemma}
\label{lemme1}
The following estimate holds on the semi-discrete approximation of\eqref{homog}. If $ v \in \mathcal{D}((-A)^{\beta/2})$
\begin{eqnarray}
 \label{form4}
 \Vert u(t)-u_{h}(t)\Vert &=&\Vert T_{h}(t) v\Vert \leq C h^{r} t^{-(r-\beta)/2}\Vert v \Vert_{\beta}\;\;\;\;\; r \in \left\lbrace 1,2 \right\rbrace\;\;\beta \leq r,
\end{eqnarray}
where $r$ is  related to \eqref{regulard}.

{\bf Proof}

The proof for $r=2$ and $\beta=0$ can be found in [\cite{lions}, Theorem 7.1, page 817].\\
Set
\begin{eqnarray}
  u_{h}(t)-u(t) = \left(u_{h}(t)-R_{h}u(t)\right)+ \left(R_{h}u(t)-u(t)\right)\equiv \theta(t) + \rho(t).
\end{eqnarray}
It is well known \cite{Stig} that $A_{h}R_{h}=P_{h}A$. Indeed for $v \in \mathcal{D}(A),\;\chi\in V_{h} $ we have 
\begin{eqnarray*}
 (P_{h}Av,\chi)&=&(A v,\chi)\qquad \qquad \qquad \qquad\;\;\; (\text{by definition of $P_{h}$})\\
                        &=& ( AR_{h}v,\chi)\qquad \qquad\qquad \qquad (\text{by definition of $R_{h}$})\\
                          &=& (A_{h}R_{h}v,\chi) \qquad \qquad \qquad \qquad (\text{since  $R_{h}v \in V_{h} $})
\end{eqnarray*}
thus $A_{h}R_{h}=P_{h}A $.
We therefore have the following equation in $\theta$
\begin{eqnarray*}
 \theta_{t}=A_{h}\theta -P_{h}D_{t}\rho.
\end{eqnarray*}
Hence
\begin{eqnarray*}
 \theta(t)=S_{h}(t)\theta(0)-\int_{0}^{t}S_{h}(t-s)P_{h} D_{s}\rho ds.
\end{eqnarray*}

Splitting the integral up into two intervals
and integration by parts over the first interval yields
\begin{eqnarray*}
 \theta(t)= S_{h}(t)\theta(0)+S_{h}(t)P_{h}\rho(0)-S_{h}(t/2)P_{h}\rho(t/2) +\int_{0}^{t/2}\left(D_{s}S_{h}(t-s) \right) P_{h}\rho(s)ds -\int_{t/2}^{t}S_{h}(t-s)P_{h}D_{s}\rho(s)ds,
\end{eqnarray*}
with $D_{s}=\partial/\partial s$.
Since $\theta (t) \in V_{h}$ we therefore have $P_{h}\theta (t)=\theta (t)$, then
\begin{eqnarray*}
\theta(t)= S_{h}(t)P_{h}T_{h}(0)v-S_{h}(t/2)P_{h}\rho(t/2)+\int_{0}^{t/2}\left(D_{s}S_{h}(t-s)\right)P_{h}\rho(s)ds-\int_{t/2}^{t}S_{h}(t-s)P_{h}D_{s}\rho(s)ds.
\end{eqnarray*}
Since $$ P_{h}T_{h}(0)v=P_{h}(v-P_{h}v)=0, $$  we therefore have
 \begin{eqnarray*}
  \theta(t)=-S_{h}(t/2)P_{h}\rho(t/2)+\int_{0}^{t/2}D_{s}S_{h}(t-s)P_{h}\rho(s)ds-\int_{t/2}^{t}S_{h}(t-s)P_{h} D_{s}\rho(s)ds.
 \end{eqnarray*}
Using the fact that  $S_{h}$ and $P_{h}$ are uniformly bounded independently of $h$ with the smoothing property of $S_{h}$ in \propref{prop1} yields
\begin{eqnarray*}
 \Vert \theta(t)\Vert \leq C \left( \Vert \rho(t/2)\Vert +  \int_{0}^{t/2} (t-s)^{-1} \Vert \rho(s)\Vert ds + \int_{t/2}^{t} \Vert D_{s}\rho(s)\Vert ds\right).
\end{eqnarray*}

The estimate \eqref{regulard} with the smoothing property of $S(t)$ in \propref{prop1} yields 
\begin{eqnarray*}
\left\lbrace \begin{array}{l}
\Vert \rho(t)\Vert \leq  C h^{r} \Vert u\Vert_{r} \leq C h^{r} t^{-(r-\beta)/2} \Vert v\Vert_{\beta}\\
\newline\\
\Vert D_{s}\rho(t)\Vert \leq  C h^{r} \Vert  D_{s}u\Vert_{r} \leq  C h^{r}t^{-1-(r-\beta)/2} \Vert v\Vert_{\beta}
,\;\;\;\ r\in \left\lbrace  1,2\right\rbrace,\;\; \beta\leq r , \;\;\; v\in \mathcal{D}((-A)^{\beta/2}). 
 \end{array}\right.
\end{eqnarray*}
Then
\begin{eqnarray*}
 \Vert \theta(t)\Vert \leq C h^{r} t^{-(r-\beta)/2}\Vert v\Vert_{\beta}  + C h^{r} \Vert v\Vert_{\beta} \left(\int_{0}^{t/2} (t-s)^{-1} s^{-(r-\beta)/2}ds+\int_{t/2}^{t} s^{-1-(r-\beta)/2}
ds\right).
\end{eqnarray*}
Since
$$\int_{0}^{t/2} (t-s)^{-1} s^{-(r-\beta)/2}ds+\int_{t/2}^{t} s^{-1-(r-\beta)/2}ds\leq  C t^{-(r-\beta)/2}, $$
we therefore have
\begin{eqnarray*}
 \Vert T_{h}(t)v \Vert \leq \Vert \theta(t)\Vert +\Vert \rho(t)\Vert \leq C h^{r} \;t^{-(r-\beta)/2}\Vert v\Vert_{\beta}.
\end{eqnarray*}



\end{lemma}
Our second preliminary lemma concerns the mild solution SPDE of (\ref{adr}).
\begin{lemma}
\label{lemme2}
Let $X$ be the  mild solution given in (\ref{eq1}).
Suppose that \assref{assumption1} and  \assref{assumption2}  hold. Let  $0 \leq \gamma < 1 $, $t_{1}, t_{2} \in [0,T]$ be so that $t_{1}<
t_{2}$. If $X_{0} \in L_{2}(\mathbb{D},\mathcal{D}((-A)^{\gamma}))$ then we have the following estimate, 
\begin{eqnarray*}
\mathbf{E}\Vert X(t_{2})- X(t_{1}) \Vert^{2} &\leq&  C
(t_{2}-t_{1})^{\gamma }\left(\mathbf{E} \Vert X_{0}\Vert_{\gamma}^{2}+\mathbf{E} \left(\underset{0\leq s\leq T}{\sup} \left(1+
 \Vert X(s)\Vert \right)\right)^{2}+ \left(\underset{0\leq s\leq t_{1}}{\sup}\mathbf{E} \left(1+\Vert X(s)\Vert\right)^{2} \right)\right).
\end{eqnarray*}
\end{lemma}
{\bf Proof}
Consider the difference
\begin{eqnarray*}
\lefteqn{  X(t_{2})- X(t_{1})} & &\\
&=&  \left(S(t_{2})-S(t_{1})\right)X_{0}+\left(
  \int_{0}^{t_{2}}S(t_{2}-s)F(X(s))ds-\int_{0}^{t_{1}}S(t_{1}-s)F(X(s))ds\right) \\  
& &  +\left(\int_{0}^{t_{2}}S(t_{2}-s)B(X)dW(s)-\int_{0}^{t_{1}}S(t_{1}-s)B(X)dW(s)\right)\\ 
&=& I +II +III, 
\end{eqnarray*}
so that 
$ \mathbf{E}  \Vert X(t_{2})- X(t_{1})\Vert^{2} \leq 3 (\mathbf{E}  \Vert I \Vert^{2} + \mathbf{E} \Vert II \Vert^{2} + \mathbf{E} \Vert III \Vert^{2}).$
We estimate each of the terms $I, II$ and $III$. Estimation of the terms $I$ and $II$ are similar to ones in [\cite{GTambue,GTambueexpo}, Lemma 3.2] with additive noise.
 Using \propref{prop1} as in \cite{GTambue,GTambueexpo} yields
 \begin{eqnarray*}
 \Vert I \Vert
&=&\Vert S(t_{1})(-A)^{-\gamma/2}(\text{I}-S(t_{2}-t_{1})) (-A)^{\gamma/2} X_{0} \Vert 
\quad \leq \quad   C (t_{2}-t_{1})^{\gamma/2} \Vert X_{0} \Vert_{\gamma}.
\end{eqnarray*}
and 
\begin{eqnarray*}
\mathbf{E} \Vert II\Vert^{2}\leq C (t_{2}-t_{1})^{2 \gamma}  \mathbf{E}\left(\underset{0\leq s\leq T}{\sup} (1+ \Vert X(s)\Vert) \right)^{2}.
\end{eqnarray*}

For term  $III$, we have 
\begin{eqnarray*}
 III 
&=&\int_{0}^{t_{1}}\left( S(t_{2}-s)-S(t_{1}-s)\right)B(X)dW(s)+ \int_{t_{1}}^{t_{2}}S(t_{2}-s)B(X)dW(s)
  \quad = \quad III_{1} +III_{2}.
\end{eqnarray*}
Using the Ito isometry property yields
\begin{eqnarray*}
\mathbf{E} \Vert III_{1} \Vert ^{2}&=& \mathbf{E}\Vert \int_{0}^{t_{1}}\left( S(t_{2}-s)-S(t_{1}-s)\right)B(X)dW(s)\Vert^{2}\\
 &=&  \int_{0}^{t_{1}}\mathbf{E} \Vert \left(S(t_{2}-s)-S(t_{1}-s)\right)B(X) Q^{1/2}\Vert_{HS}^{2} ds.
\end{eqnarray*}
 Using \assref{assumption2} and \propref{prop1} yields 
\begin{eqnarray*}
 \mathbf{E} \Vert III_{1} \Vert ^{2}&\leq & C \left( \int_{0}^{t_{1}} \Vert ( S(t_{2}-s)-S(t_{1}-s))\Vert_{L(L^{2}(\Omega))}^{2} ds \right) \left(\underset{0\leq s\leq t_{1}}{\sup}\mathbf{E} \left(1+\Vert X(s)\Vert\right)^{2} \right) \\ 
&=  & C\left(\int_{0}^{t_{1}} \Vert S(t_{1}-s)(-A)^{\gamma/2} (-A)^{-\gamma/2} ( \textbf{I}-S(t_{2}-t_{1}))\Vert_{L(L^{2}(\Omega))}^{2} ds\right)  \left(\underset{0\leq s\leq t_{1}}{\sup} \mathbf{E} \left(1+\Vert X(s)\Vert\right)^{2} \right)\\
&=& C \left( \int_{0}^{t_{1}} \Vert (-A)^{\gamma/2}  S(t_{1}-s) (-A)^{-\gamma/2} ( \textbf{I}-S(t_{2}-t_{1}))\Vert_{L(L^{2}(\Omega))}^{2} ds \right) \left(\underset{0\leq s\leq t_{1}}{\sup}\mathbf{E}\left(1+\Vert X(s)\Vert\right)^{2} \right) \\
&\leq & C (t_{2}-t_{1})^{\gamma} \left( \int_{0}^{t_{1}} (t_{1}-s)^{-\gamma}ds \right)  \left(\underset{0\leq s\leq t_{1}}{\sup}\mathbf{E} \left(1+\Vert X(s)\Vert\right)^{2} \right) \\
&\leq& C (t_{2}-t_{1})^{\gamma} \left(\underset{0\leq s\leq t_{1}}{\sup}\mathbf{E} \left(1+\Vert X(s)\Vert\right)^{2}\right) .
\end{eqnarray*}
Let us estimate  $\mathbf{E} \Vert III_{2} \Vert$. The Ito isometry
again,  with  the boundedness of $S$ and \assref{assumption2} yields
\begin{eqnarray*}
 \mathbf{E} \Vert III_{2} \Vert ^{2}&=& \mathbf{E}\Vert \int_{t_{1}}^{t_{2}}S(t_{2}-s))B(X)dW(s)\Vert^{2}\\
&=& \int_{t_{1}}^{t_{2}}\mathbf{E}  \Vert S(t_{2}-s))B(X(s))\Vert_{L_{0}^{2}}^{2}ds\\
& =& \;\int_{t_{1}}^{t_{2}} \Vert S(t_{2}-s)\Vert_{L(L^{2}(\Omega))} \, ds \left(\underset{0\leq s\leq t_{1}}{\sup}\mathbf{E}\left(1+\Vert X(s)\Vert\right)^{2}  \right)  \\
 &\leq & \; C (t_{2}-t_{1})\left(\underset{0\leq s\leq t_{1}}{\sup}\mathbf{E} \left(1+\Vert X(s)\Vert\right)^{2}  \right).
\end{eqnarray*}
Hence 
\begin{eqnarray*}
\mathbf{E} \Vert III\Vert^{2}\leq 2( \mathbf{E} \Vert III_{1} \Vert^{2} +\mathbf{E} \Vert III_{2} \Vert^{2} )\leq C (t_{2}-t_{1})^{\gamma}. 
\end{eqnarray*}
Combining our estimates of  $\mathbf{E} \Vert I\Vert^{2},\mathbf{E}
\Vert II\Vert^{2}$ and $\mathbf{E} \Vert III\Vert^{2}$ ends the  lemma.

\subsection{Proof of \thmref{th1} for the scheme SETDM1}
{\bf Proof}
Set 
\begin{eqnarray*}
 X(t_{m})&=&S(t_{m})X_{0}+\underset{k=0}{\sum^{m-1}}\int_{t_{k}}^{t_{k+1}}S(t_{m}-s)F(X(s))ds+\int_{0}^{t_{m}} S(t_{m}-s)B(X(t_{m}))d W(s)\\
          &=& \overline{X}(t_{m})+ O(t_{m}).
\end{eqnarray*}
%
Recall that
\begin{eqnarray*}
 X_{m}^{h}&=& e^{\Delta t A_{h}}X_{m-1}^{h}+A_{h}^{-1}\left( e^{\Delta t A_{h}}-I\right)P_{h}F(X_{m-1}^{h})+ \int_{t_{m-1}}^{t_{m}} e^{(t_{m}-s)A_{h}}P_{h}B(X_{m-1}^{h})d W(s)\\
 &=& e^{\Delta t A_{h}}X_{m-1}^{h}+\int_{0}^{\Delta t} e^{(\Delta t- s)A_{h}}P_{h}F(X_{m-1}^{h})ds + \int_{t_{m-1}}^{t_{m}} e^{(t_{m}-s)A_{h}}P_{h}B(X_{m-1}^{h})d W(s)\\
&=& S_{h}(t_{m}) P_{h}X_{0}+\underset{k=0}{\sum^{m-1}}\left(\int_{t_{k}}^{t_{k+1}} S_{h}(t_{m}-s)P_{h}F(X_{k}^{h})ds +\int_{t_{k}}^{t_{k+1}} S_{h}(t_{m}-s)P_{h}B(X_{k}^{h})dW(s) \right)\\
&=& S_{h}(t_{m}) P_{h}X_{0}+\underset{k=0}{\sum^{m-1}}\left(\int_{t_{k}}^{t_{k+1}} S_{h}(t_{m}-s)P_{h}F(X_{k}^{h})ds\right)+ \underset{k=0}{\sum^{m-1}} \int_{t_{k}}^{t_{k+1}} S_{h}(t_{m}-s)P_{h}B(X_{k}^{h})d W(s)\\
&=& Z_{m}^{h} +O_{m}^{h},
\end{eqnarray*}
with
\begin{eqnarray*}
 Z_{m}^{h}&=& S_{h}(t_{m}) P_{h}X_{0}+\underset{k=0}{\sum^{m-1}}\left(\int_{t_{k}}^{t_{k+1}} S_{h}(t_{m}-s)P_{h}F(X_{k}^{h})ds\right).
\end{eqnarray*}
We examine the error 
 \begin{eqnarray}
  X(t_{m})-X_{m}^{h}
  &=& \overline{X}(t_{m})+ O(t_{m})-X_{m}^{h}\nonumber \\
  &=&\overline{X}(t_{m}) + O(t_{m})-\left(Z_{m}^{h}+O_{m}^{h}\right)\nonumber \\ 
  &=& \left(\overline{X}(t_{m}) -Z_{m}^{h}\right)+\left( O(t_{m})-O_{m}^{h}\right)\nonumber \\
  &=& I +II \label{eq:IIIIII},
\end{eqnarray}
thus
\begin{eqnarray}
  \mathbf{E} \Vert X(t_{m})-X_{m}^{h}\Vert^{2}  \leq 2\left(\mathbf{E} \Vert I\Vert^{2} + \mathbf{E} \Vert II\Vert^{2}\right).
\end{eqnarray}

We follow the approach in \cite{GTambueexpo}. Let us estimate  the
first term $\textbf{E}\Vert I\Vert^{2}$.  
Using the definition of $T_h$ from \eqref{form1}, the first
term $I$ can be expanded 
\begin{eqnarray}
I &=&  T_{h}X_{0} + \underset{k=0}{\sum^{m-1}} \int_{t_{k}
  t}^{t_{k+1}} S( t_{m}-s) F(X(s))-S_{h}(t_{m}-s)
P_{h}F(X_{k}^{h})ds \nonumber \\  
 &=& T_{h}X_{0} + \underset{k=0}{\sum^{m-1}} \int_{t_{k}}^{t_{k+1}}S_{h}(t_{m}-s) P_{h}(F(X(t_{k}))
-F(X_{k}^{h})))ds\nonumber \\ 
&& +\underset{k=0}{\sum^{m-1}} \int_{t_{k}}^{t_{k+1}}
S_{h}(t_{m}-s)P_{h}(F(X(s))-F (X(t_{k})))ds\nonumber \\ 
&& +\underset{k=0}{\sum^{m-1}} \int_{t_{k}}^{t_{k+1}}(S(t_{m}-s)-S_{h}(t_{m}-s)P_{h})F(X(s))ds\nonumber \\ 
&=& I_{1}+I_{2}+I_{3}+I_{4}.
\label{eq:I1toI5}
 \end{eqnarray}

Then
\begin{eqnarray*}
\mathbf{E} \Vert I\Vert^{2}\leq  4 \left(\mathbf{E} \Vert I_{1} \Vert^{2}+\mathbf{E} \Vert I_{2} \Vert^{2}+\mathbf{E} \Vert I_{3} \Vert^{2}+\mathbf{E} \Vert I_{4} \Vert^{2}\right).
\end{eqnarray*}
Let us estimate $I_1$, for $0 \leq \gamma<1$ with $ 2\gamma \leq r $
and $r \in \left\lbrace 1,2\right\rbrace $,   if $X_{0} \in
L_{2}(\mathbb{D},\mathcal{D}((-A)^{\gamma}))$, Lemma \ref{lemme1} with
$\beta=2 \gamma$  yields 
\begin{eqnarray*}
 \mathbf{E} \Vert I_{1} \Vert^{2}\leq C 
 t_{m}^{-(r-2\gamma)}h^{2 r} \left(\textbf{E} \Vert
   X_{0}\Vert_{2\gamma}^{2}\right),
\end{eqnarray*}
and if  $X_{0} \in \mathcal{D}(-A)$ we have 
\begin{eqnarray*}
 \mathbf{E} \Vert I_{1} \Vert^{2}\leq C h^{4} \textbf{E} \Vert
   X_{0}\Vert_{2}^{2}.
\end{eqnarray*}

For $I_{2}$, using \assref{assumption1}, triangle inequality as well as 
the fact that $S_{h}(t) $ and $P_{h}$ are  bounded operators with
Fubini's theorem yields 
\begin{eqnarray*}
\mathbf{E} \Vert I_{2} \Vert^{2} &\leq & C
 m \underset{k=0}{\sum^{m-1}}\mathbf{E} \Vert \int_{t_{k}}^{t_{k+1}}
 S_{h}(t_{m}-s)P_{h}\left(  F(X(t_{k}))-F(X_{k}^{h})\right) ds\Vert^{2} \\ 
&\leq & C m \underset{k=0}{\sum^{m-1}}\mathbf{E}\left(   \int_{t_{k}}^{t_{k+1}}
 \Vert F(X(t_{k}))-F(X_{k}^{h})\Vert ds\right)^{2} \\ 
 &\leq & C m \,\Delta t \,\underset{k=0}{\sum^{m-1}}\int_{t_{k}}^{t_{k+1}
  }\left(\mathbf{E}\Vert X(t_{k}
   )-X_{k}^{h}\Vert^{2}\right)ds\\
&\leq & C \underset{k=0}{\sum^{m-1}}\int_{t_{k}}^{t_{k+1}
  }\left(\mathbf{E}\Vert X(t_{k}
   )-X_{k}^{h}\Vert^{2}\right)ds.
\end{eqnarray*}

Once again using the Lipschitz  condition, triangle inequality, the fact that  $S_{h}$ and $P_{h}$ are bounded but with Lemma \ref{lemme2} yields 
\begin{eqnarray*}
 (\mathbf{E} \Vert I_{3} \Vert^{2})^{1/2} &\leq &\underset{k=0}{\sum^{m-1}}\int_{t_{k}}^{t_{k+1}} (\textbf{E}  \Vert S_{h}(t_{m}-s)P_{h}(F(X(s))-F(X(t_{k}))\Vert^{2})^{1/2}ds\\
&\leq & C \underset{k=0}{\sum^{m-1}}\int_{t_{k}}^{t_{k+1}}(\textbf{E}  \Vert F(X(s))-F(X(t_{k}))\Vert)^{1/2}ds\\
&\leq&  C \left (\underset{k=0}{\sum^{m-1}} \int_{t_{k}}^{t_{k+1}}(s- t_{k})^{\gamma/2}ds\right)\;\\
&& \times \left(\mathbf{E} \Vert X_{0}\Vert_{\gamma}^{2}+\textbf{E} \left(\underset{0\leq s\leq T}{\sup} (1+ \Vert X(s)\Vert)\right)^{2}+ \left(\underset{0\leq s\leq T}{\sup}\textbf{E}\left( 1+ \Vert X(s)\Vert \right)^{2} \right)\right)^{1/2} \\
                                   &\leq&  C  \Delta t^{\gamma/2} \left(\textbf{E} \Vert X_{0}\Vert_{\gamma}^{2}+\textbf{E} \left(\underset{0\leq s\leq T}{\sup} (1+ \Vert X(s)\Vert)\right)^{2}+ \left(\underset{0\leq s\leq T}{\sup}\textbf{E} \left( 1+ \Vert X(s)\Vert \right)^{2} \right)\right)^{1/2},
\end{eqnarray*}
 thus
\begin{eqnarray*}
 \mathbf{E}\Vert I_{3} \Vert^{2} \leq C \Delta t^{\gamma} .
\end{eqnarray*}

If  $X_{0} \in L_{2}\left( \mathbb{D}, \mathcal{D}(-A)\right)$  we obviously have 
$\mathbf{E} \Vert I_{3} \Vert^{2} \leq C (\Delta t)^{1-\epsilon}$
by taking $\gamma = 1-\epsilon$ in Lemma  \ref{lemme2},\;$\epsilon \in
(0,1/2)$ small enough. 
Let us estimate $\left(\textbf{E} \Vert I_{4} \Vert^{2}\right)^{1/2}$. For $r=1, \beta=0$, using  \lemref{lemme1} yields
\begin{eqnarray*}
 \left(\mathbf{E} \Vert I_{4} \Vert^{2}\right)^{1/2} &\leq&
 \underset{k=0}{\sum^{m-1}}\int_{t_{k}}^{t_{k+1}}  \left(\textbf{E}
   \Vert T_{h}(t_{m}-s) F(X(s))\Vert^2 \right)^{1/2} ds \\
& \leq &  C h \underset{0\leq s\leq T}{\sup}\left( \textbf E \Vert F( X(s)) \Vert^{2}\right)^{1/2}
\left(\int_{0}^{t_{m}}\left(t_{m}-s \right)^{-1/2} \right)\\
& \leq &  C h,
\end{eqnarray*}
thus
\begin{eqnarray*}
 \mathbf{E} \Vert I_{4} \Vert^{2}\leq C h^{2}.
\end{eqnarray*}
If $  F(L_{2}(\mathbb{D},\mathcal{D}((-A)^{\alpha'})))\subset
L_{2}(\mathbb{D},\mathcal{D}((-A)^{\alpha'})),\; \alpha' \in
(0,\gamma/10)$ small enough, for $r=2, \beta=2 \alpha'$, using
\lemref{lemme1} yields 
\begin{eqnarray*}
 \left(\mathbf{E} \Vert I_{4} \Vert^{2}\right)^{1/2} &\leq&
 \underset{k=0}{\sum^{m-1}}\int_{t_{k}}^{t_{k+1}}  \left(\textbf{E}
   \Vert T_{h}(t_{m}-s) F(X(s))\Vert_{\beta}^2 \right)^{1/2} ds \\
& \leq &  C h^{2} \underset{0\leq s\leq T}{\sup}\left( \mathbf E \Vert F( X(s)) \Vert_{\beta}^{2}\right)^{1/2}
\left(\int_{0}^{t_{m}}\left(t_{m}-s \right)^{-1+\beta/2} ds\right)\\
& \leq &  C h^{2},
\end{eqnarray*}
thus
\begin{eqnarray*}
 \mathbf{E} \Vert I_{4} \Vert^{2}\leq C h^{4}.
\end{eqnarray*}
Combining the previous estimates yields:
For  $X_{0}\in L_{2}\left(\mathbb{D}, \mathcal{D}((-A)^{\gamma})\right),\;\; 1/2\leq \gamma <1$ 
\begin{eqnarray*}
 \mathbf{E} \Vert I \Vert^{2} \leq C \left(h^{2} +\Delta t^{\gamma} + \underset{k=0}{\sum^{m-1}}\int_{t_{k}}^{t_{k+1}
  }\left(\mathbf{E}\Vert X(t_{k}
   )-X_{k}^{h}\Vert^{2}\right)ds\right).\\
\end{eqnarray*}
For  $X_{0}\in L_{2}\left(\mathbb{D}, \mathcal{D}((-A)^{\gamma})\right),\; 0 < \gamma \leq 1/2$ 
\begin{eqnarray*}
 \mathbf{E} \Vert I \Vert^{2} \leq C \left(t_{m}^{-1+2\gamma} h^{2} +\Delta t^{\gamma} + \underset{k=0}{\sum^{m-1}}\int_{t_{k}}^{t_{k+1}
  }\left(\mathbf{E}\Vert X(t_{k}
   )-X_{k}^{h}\Vert^{2}\right)ds\right).\\
\end{eqnarray*}
For  $X_{0}\in
L_{2}\left(\mathbb{D},\mathcal{D}((-A)^{\gamma})\right),\; 0\leq
\gamma < 1$ and  if  $
F(L_{2}(\mathbb{D},\mathcal{D}((-A)^{\alpha'})))\subset
L_{2}(\mathbb{D},\mathcal{D}((-A)^{\alpha'})),\; \alpha' \in
(0,\gamma/10)$ small enough 
 \begin{eqnarray*}
 \mathbf{E} \Vert I \Vert^{2} \leq C \left( t_{m}^{-2+2\gamma} h^{4} +\Delta t^{\gamma} + \underset{k=0}{\sum^{m-1}}\int_{t_{k}}^{t_{k+1}
  }\left(\mathbf{E}\Vert X(t_{k}
   )-X_{k}^{h}\Vert^{2}\right)ds\right).\\
\end{eqnarray*}
For $X_{0}\in L_{2}\left(\mathbb{D}, \mathcal{D}(-A)\right) $ and   if
$  F(L_{2}(\mathbb{D},\mathcal{D}((-A)^{\alpha'})))\subset
L_{2}(\mathbb{D},\mathcal{D}((-A)^{\alpha'})),\; \alpha' \in
(0,\gamma/10)$ small enough, 
\begin{eqnarray*}
 \mathbf{E} \Vert I \Vert^{2} \leq C \left( h^{4} +\Delta
 t^{1-\epsilon} + \underset{k=0}{\sum^{m-1}}\int_{t_{k}}^{t_{k+1}
 }\left(\mathbf{E}\Vert X(t_{k} 
   )-X_{k}^{h}\Vert^{2}\right)ds\right),
\end{eqnarray*}
with  $\epsilon \in (0,1/2)$ small  enough.

Let us estimate $\mathbf{E} \Vert II \Vert^{2}$, we follow the same approach as in \cite{Yn:05}. Note that in the case of additive noise the 
estimation is straightforward and smooth  noise improve the accuracy (see \cite{KLNS} and \figref{FIG0022} in \secref{simulation}).
For multiplicative noise we have 
\begin{eqnarray}
II &=&  \int_{0}^{t_{m}} S(t_{m}-s)B(X(s))d W(s)- \underset{k=0}{\sum^{m-1}} \int_{t_{k}}^{t_{k+1}} S_{h}(t_{m}-t_{k})P_{h}B(X_{k}^{h})d W(s)\nonumber \\
 &=& \underset{k=0}{\sum^{m-1}} \int_{t_{k}}^{t_{k+1}}S_{h}(t_{m}-t_{k}) P_{h}\left(B(X(t_{k}))
-B(X_{k}^{h}))\right) d W(s)\nonumber \\ 
&& +\underset{k=0}{\sum^{m-1}} \int_{t_{k}}^{t_{k+1}}
S_{h}(t_{m}-t_{k})P_{h}(B(X(s))-B(X(t_{k})))d W(s)\nonumber \\ 
&& +\underset{k=0}{\sum^{m-1}} \int_{t_{k}}^{t_{k+1}}(S(t_{m}-t_{k})-S_{h}(t_{m}-t_{k})P_{h})B(X(s))d W(s)\nonumber \\ 
&& +\underset{k=0}{\sum^{m-1}} \int_{t_{k}}^{t_{k+1}} \left(S(t_{m} 
 -s )-S(t_{m} -t_{k}))B(X(s))\right)dW(s)\nonumber \\ 
&=& II_{1}+II_{2}+II_{3}+II_{4}.
\label{eq:I1eto}
 \end{eqnarray}
Then
\begin{eqnarray*}
\mathbf{E} \Vert II\Vert^{2}\leq 4 \left(\mathbf{E} \Vert II_{1} \Vert^{2}+\mathbf{E} \Vert II_{2} \Vert^{2}+\mathbf{E} \Vert II_{3} \Vert^{2}+\mathbf{E} \Vert II_{4} \Vert^{2}\right).
\end{eqnarray*}
Let us estimate each term. Using the Ito isometry, the boundedness of
$S_{h}$  and $P_{h}$, and \assref{assumption2} yields 
\begin{eqnarray*}
\mathbf{E} \Vert II_{1} \Vert^{2} &=& \mathbf{E}\Vert \underset{k=0}{\sum^{m-1}} \int_{t_{k}}^{t_{k+1}}S_{h}(t_{m}-t_{k}) P_{h}\left(B(X(t_{k}))
-B(X_{k}^{h})\right) d W(s)\Vert^{2}\\
 &=& \underset{k=0}{\sum^{m-1}} \mathbf{E}\Vert \int_{t_{k}}^{t_{k+1}}S_{h}(t_{m}-t_{k}) P_{h}\left(B(X(t_{k}))
-B(X_{k}^{h})\right) d W(s)\Vert^{2}\\
&\leq & C\underset{k=0}{\sum^{m-1}} \int_{t_{k}}^{t_{k+1}} \mathbf{E}\Vert B(X(t_{k}))
-B(X_{k}^{h})\Vert_{L_{0}^{2}} ds\\
&\leq& C \underset{k=0}{\sum^{m-1}} \int_{t_{k}}^{t_{k+1}}\mathbf{E}\Vert X(t_{k})
-X_{k}^{h}\Vert^{2} ds.
\end{eqnarray*}
For $X_{0}\in L_{2}\left(\mathbb{D}, \mathcal{D}\left( (-A)^{\gamma}\right)\right)$, using \lemref{lemme2} and \assref{assumption2} yields
\begin{eqnarray*}
 \mathbf{E} \Vert II_{2} \Vert^{2} &=& \mathbf{E}\Vert \underset{k=0}{\sum^{m-1}} \int_{t_{k}}^{t_{k+1}}
S_{h}(t_{m}-t_{k})P_{h}(B(X(s))-B(X(t_{k})))d W(s)\Vert^{2}\\
 &=& \underset{k=0}{\sum^{m-1}} \int_{t_{k}}^{t_{k+1}}\mathbf{E}\Vert
S_{h}(t_{m}-t_{k})P_{h}(B(X(s))-B(X(t_{k})))\Vert_{L_{2}^{0}}^{2}ds\\
&\leq& C \underset{k=0}{\sum^{m-1}}\int_{t_{k}}^{t_{k+1}} \mathbf{E}\Vert X(s)
-X(t_{k})\Vert^{2}ds \\
&\leq& C \left(\underset{k=0}{\sum^{m-1}}\int_{t_{k}}^{t_{k+1}} (s-t_{k})^{\gamma}ds \right)\\
&& \times\left(\mathbf{E} \Vert X_{0}\Vert_{\gamma}^{2}+\mathbf{E} \left(\underset{0\leq s\leq T}{\sup} (1 + \Vert X(s)\Vert)\right)^{2}+ \left(\underset{0\leq s\leq T}{\sup}\mathbf{E} \Vert\left(1+ X(s)\Vert\right)^{2} \right)\right)\\
&\leq& C \Delta t^{\gamma}. 
\end{eqnarray*}

For   $X_{0}\in L_{2}\left(\mathbb{D}, \mathcal{D}(-A)\right)$, taking
$\gamma=1-\epsilon$, with $\epsilon$ small enough yields 
\begin{eqnarray*}
 \mathbf{E} \Vert II_{2} \Vert^{2} \leq  C \Delta t^{1-\epsilon}. 
\end{eqnarray*}
Let us estimate $\mathbf{E} \Vert II_{3} \Vert^{2}$. By Ito's isometry and \lemref{lemme1}, we have
\begin{eqnarray*}
 \mathbf{E} \Vert II_{3} \Vert^{2}&=&  \mathbf{E}\Vert\underset{k=0}{\sum^{m-1}} \int_{t_{k}}^{t_{k+1}}(S(t_{m}-t_{k})-S_{h}(t_{m}-t_{k})P_{h})B(X(s))d W(s)\Vert^{2}\\
 &=& \underset{k=0}{\sum^{m-1}} \int_{t_{k}}^{t_{k+1}}\mathbf{E}\Vert(S(t_{m}-t_{k})-S_{h}(t_{m}-t_{k})P_{h})B(X(s))\Vert_{L_{0}^{2}}^{2}ds\\
&=& \underset{k=0}{\sum^{m-1}} \int_{t_{k}}^{t_{k+1}}\mathbf{E} \Vert T_{h}(t_{m}-t_{k})B(X(s))\Vert_{L_{0}^{2}}^{2}ds.\\
\end{eqnarray*}
Indeed  using \lemref{lemme1} and \assref{assumption2}, for $
b(L_{2}(\mathbb{D},\mathcal{D}((-A)^{\alpha})))\subset
L_{2}(\mathbb{D},\mathcal{D}((-A)^{\alpha})),\; \alpha \in
(0,\gamma/10)$ small enough,  we have 
\begin{eqnarray*}
 \Vert T_{h}(t_{m}-t_{k})B(X(s))\Vert_{L_{0}^{2}}^{2}&=&\underset{i \in \mathbb{N}}{\sum^{\infty}} \Vert T_{h}(t_{m}-t_{k})b(X(s))Q^{1/2}e_{i}\Vert^{2} \\
                                              &\leq & \underset{i \in \mathbb{N}}{\sum^{\infty}} \Vert T_{h}(t_{m}-t_{k})b(X(s))\Vert^{2} \Vert Q^{1/2}e_{i}\Vert^{2} \\
                                              &\leq & C  h^{2} (t_{m}-t_{k})^{-1+\alpha}\Vert b(X(s))\Vert_{\alpha}^{2} \mathbf{Tr}(Q),
\end{eqnarray*}
thus 
\begin{eqnarray*}
 \mathbf{E} \Vert II_{3} \Vert^{2} \leq  C  h^{2} \mathbf{Tr}(Q)
 \underset{0\leq s\leq T}{\sup}\mathbf{E} \Vert
 b(X(s))\Vert_{\alpha}^{2} \left( \Delta t^{\alpha}\underset{
   k=0}{\sum^{m-1}} (m-k)^{-1+\alpha}\right), 
\end{eqnarray*}
since  
$$
 \underset{ k=0}{\sum^{m-1}} (m-k)^{-1+\alpha},
$$
is the discrete form of $$ \Delta t^{\alpha}\int_{0}^{m-1}(m-s)^{-1+\alpha}ds  \leq  \Delta t^{\alpha}M^{\alpha}=T^{\alpha},$$
we therefore have 
\begin{eqnarray*}
 \mathbf{E} \Vert II_{3} \Vert^{2} \leq  C  h^{2} \mathbf{Tr}(Q) \underset{0\leq s\leq T}{\sup}\mathbf{E} \Vert b(X(s))\Vert_{\alpha}^{2}. 
\end{eqnarray*}
For $  b(L_{2}(\mathbb{D},\mathcal{D}((-A))))\subset
L_{2}(\mathbb{D},\mathcal{D}((-A)))$, we obviously have using
\lemref{lemme1}
\begin{eqnarray*}
 \mathbf{E} \Vert II_{3} \Vert^{2} \leq  C  h^{4} \mathbf{Tr}(Q)\underset{0\leq s\leq T}{\sup}\mathbf{E} \Vert b(X(s))\Vert_{2}^{2}. 
\end{eqnarray*}
Let us estimate $\mathbf{E} \Vert II_{4} \Vert^{2}$, 
by our assumption on $b$ the following estimation holds 
\begin{eqnarray*}
 \mathbf{E} \Vert II_{4} \Vert^{2} &=& \underset{k=0}{\sum^{m-1}} \int_{t_{k}}^{t_{k+1}} \mathbf{E}\Vert (S(t_{m} 
 -s )-S(t_{m} -t_{k}))B(X(s))\Vert_{L^{2}_{0}}^{2}ds\\
&\leq & \mathbf{Tr}(Q)  \left(\underset{k=0}{\sum^{m-1}}
 \int_{t_{k}}^{t_{k+1}}\Vert (-A)^{-\alpha/2}(S(t_{m} -s )-S(t_{m}
 -t_{k}))\Vert_{L(L^{2}(\Omega))}^{2}ds\right) \left( \underset{0\leq
   s\leq T}{\sup}\mathbf{E} \Vert b(X(s))\Vert_{\alpha}^{2}\right), 
\end{eqnarray*}
since 
\begin{eqnarray*}
 \Vert (-A)^{-\alpha/2}(S(t_{m} -s )-S(t_{m} -t_{k})) 
 \Vert_{L(L^{2}(\Omega))}^{2}&=&\Vert(-A)^{(1-\alpha)/2}S(t_{m} -s )(-A)^{(-1/2)}(
 \mathbf{I}-S(s-t_{k}))\Vert_{L(L^{2}(\Omega))}^{2}\\ 
 &\leq & C (s-t_{k})(t_{m}-s)^{(\alpha-1)},
\end{eqnarray*}
thus
\begin{eqnarray*}
 \mathbf{E} \Vert II_{4} \Vert^{2}&\leq &  C
 \left(\underset{k=0}{\sum^{m-1}} \int_{t_{k}}^{t_{k+1}}
 (s-t_{k})(t_{m}-s)^{(\alpha-1)}ds \right)\left( \underset{0\leq s\leq
   T}{\sup}\mathbf{E} \Vert b(X(s))\Vert_{\alpha}^{2}\right)\\ 
 &\leq &  \Delta t  \left( \underset{0\leq s\leq T}{\sup}\mathbf{E} \Vert
 b(X(s))\Vert_{\alpha}^{2}\right)). 
\end{eqnarray*}
Combining the estimates related to $II$ yields the following. 

That for  $X_{0} \in
L_{2}\left(\mathbb{D},\mathbf{D}(-A)^{\gamma}\right)$ and  $
b(L_{2}(\mathbb{D},\mathcal{D}((-A)^{\alpha})))\subset
L_{2}(\mathbb{D},\mathcal{D}((-A)^{\alpha})),\; \alpha \in
(0,\gamma/10)$ small enough, 
\begin{eqnarray*}
\mathbf{E} \Vert II\Vert^{2}\leq   C \left(h^{2}+\Delta t^{\gamma} +\underset{k=0}{\sum^{m-1}} \int_{t_{k}}^{t_{k+1}}\mathbf{E}\Vert X(t_{k})
-X_{k}^{h}\Vert^{2} ds\right).
\end{eqnarray*}
For  $X_{0} \in L_{2}\left(\mathbb{D},\mathbf{D}(-A)^{\gamma}\right)$
and $b(L_{2}\left(\mathbb{D},\mathcal{D}(-A)\right)) \subset
L_{2}\left(\mathbb{D},\mathcal{D}(-A)\right)$ 
\begin{eqnarray*}
\mathbf{E} \Vert II\Vert^{2}\leq   C \left(h^{4}+\Delta t^{\gamma}
+\underset{k=0}{\sum^{m-1}} \int_{t_{k}}^{t_{k+1}}\mathbf{E}\Vert
X(t_{k}) -X_{k}^{h}\Vert^{2} ds\right).
\end{eqnarray*}
Combining the estimates of $ \mathbf{E} \Vert I\Vert^{2}$ and $ \mathbf{E} \Vert II\Vert^{2}$ and applying the discrete Gronwall lemma ends  the proof.
\subsection{Proof of \thmref{th1} for the scheme SETDM0}
 We just give a sketch of the mains steps.
Recall that
\begin{eqnarray*}
 Y_{m}^{h}&=& e^{\Delta t A_{h}}\left(Y_{m-1}^{h}+ \Delta t P_{h}F(Y_{m-1}^{h})\right)+ \int_{t_{m-1}}^{t_{m}} e^{\Delta t A_{h}} P_{h}B(Y_{m-1}^{h})d W(s)\\
 &=& e^{\Delta t A_{h}}Y_{m-1}^{h}+\int_{0}^{\Delta t} e^{(\Delta t)A_{h}}P_{h}F(Y_{m-1}^{h})ds + \int_{t_{m-1}}^{t_{m}} e^{\Delta A_{h}} P_{h}B(Y_{m-1}^{h})d W(s)\\
&=& S_{h}(t_{m}) P_{h}X_{0}+\underset{k=0}{\sum^{m-1}}\left(\int_{t_{k}}^{t_{k+1}} S_{h}(t_{m}-t_{k})P_{h}F(Y_{k}^{h})ds +\int_{t_{k}}^{t_{k+1}} S_{h}(t_{m}-t_{k}) P_{h}B(Y_{k}^{h})dW(s) \right)\\
&=& S_{h}(t_{m}) P_{h}X_{0}+\underset{k=0}{\sum^{m-1}}\left(\int_{t_{k}}^{t_{k+1}} S_{h}(t_{m}-t_{k})P_{h}F(Y_{k}^{h})ds\right)+ \underset{k=0}{\sum^{m-1}} \int_{t_{k}}^{t_{k+1}} S_{h}(t_{m}-t_{k})P_{h}B(Y_{k}^{h})d W(s)\\
&=& z_{m}^{h} +o_{m}^{h}.
\end{eqnarray*}
We can therefore put the estimation of the error in the form of equation \ref{eq:IIIIII}. The estimate of  
the corresponding $\mathbf{E} \Vert I\Vert^{2}$ is the  same as in Theorem \ref{th1} 
with the extra term 

$$I_{5}= \underset{k=0}{\sum^{m-1}} \int_{t_{k}}^{t_{k+1}} (S(t_{m} 
-s )-S(t_{m} -t_{k}))F(X(s))ds.$$ 
This is estimated in \cite{GTambue} Theorem 2.6 as 
$$\left(\mathbf{E}\Vert I_{5}\Vert^{2}\right)^{1/2}\leq C \left(\Delta t +\Delta t \vert \ln(\Delta t) \vert \right) \leq C \Delta t ^{\gamma/2}.$$
The estimation of $\mathbf{E} \Vert II\Vert^{2}$ is the  same as for the scheme SETDM0.

\section{Simulations}
\label{simulation}
Efficient implementation of $\varphi_{i},\; i=0,1$ can be achieved by either  the real fast \Leja points technique in~\cite{LE2,LE1,LE,GTambueexpo} or the Krylov subspace technique 
in ~\cite{kry,SID,GTambueexpo}. 
In the first example we apply the scheme to linear problem where we can construct the exact solution
for the truncated noise. The finite element method is used for space
discretization. In this example we use the real  fast L\'{e}ja point
technique to compute the exponential functions $\varphi_{i},\;
i=0,1$. We use noise with exponential correlation (see below) which is
obviously a trace class noise. 
In the second example we apply the scheme to nonlinear stochastic flow
with multiplicative noise  
in a heterogeneous media. To deal with high P\'{e}clet number flow, we
use the finite volume method for the space discretization.     
In this case we use the Krylov subspace technique  to compute the exponential
functions $\varphi_{i},\; i=0,1$, implemented in the matlab functions 
{\tt expv.m} and {\tt phiv.m} of the package Expokit~\cite{SID}. 
We compute the exponential matrix functions $\varphi_{i}$ with the Krylov
subspace technique with dimension $m=6$ and the absolute tolerance $10^{-6}$. 


In the legends of our graphs, ``SETDM1'' denotes results from the
SETDM1  scheme, ``SETDM0 '' denotes results from the SETDM0 scheme
with ``Implicit'' denotes results from  the standard  semi-implicit
Euler-Maruyama scheme. 



As a simple example consider the reaction diffusion equation  with additive noise in the time interval $[0,T]$ with diffusion coefficient $ D>0$ 
\begin{eqnarray*}
 dX=(D \varDelta X -0.5 X)dt+ dW\;\;\;\;\;\;\;\;\;\;\;\;\;\; X(0)=X_{0},\;\;\;\;\;\; \Omega=[0,L_{1}]\times [0,L_{2}]\;\
\end{eqnarray*}
with homogeneous Neumann boundary condition.
As the exact solution is known for comparison, we take  $f$   in the equation \eqref{sadr} to be linear here
\begin{eqnarray}
 f(u)=-  0.5 u.
\end{eqnarray}
The corresponding  Nemytskii operator $F$ is obtained from \eqref{nemform}.
 Of course, in general, $F$ will be nonlinear. $F$ verifies obviously \assref{assumption1}.
Here $b(x,u)=1,\; x\in \Omega,\; u \in \mathbb{R}$.

\begin{figure}[!th]
  \subfigure[]{
    \label{FIG0022a}
    \includegraphics[width=0.48\textwidth]{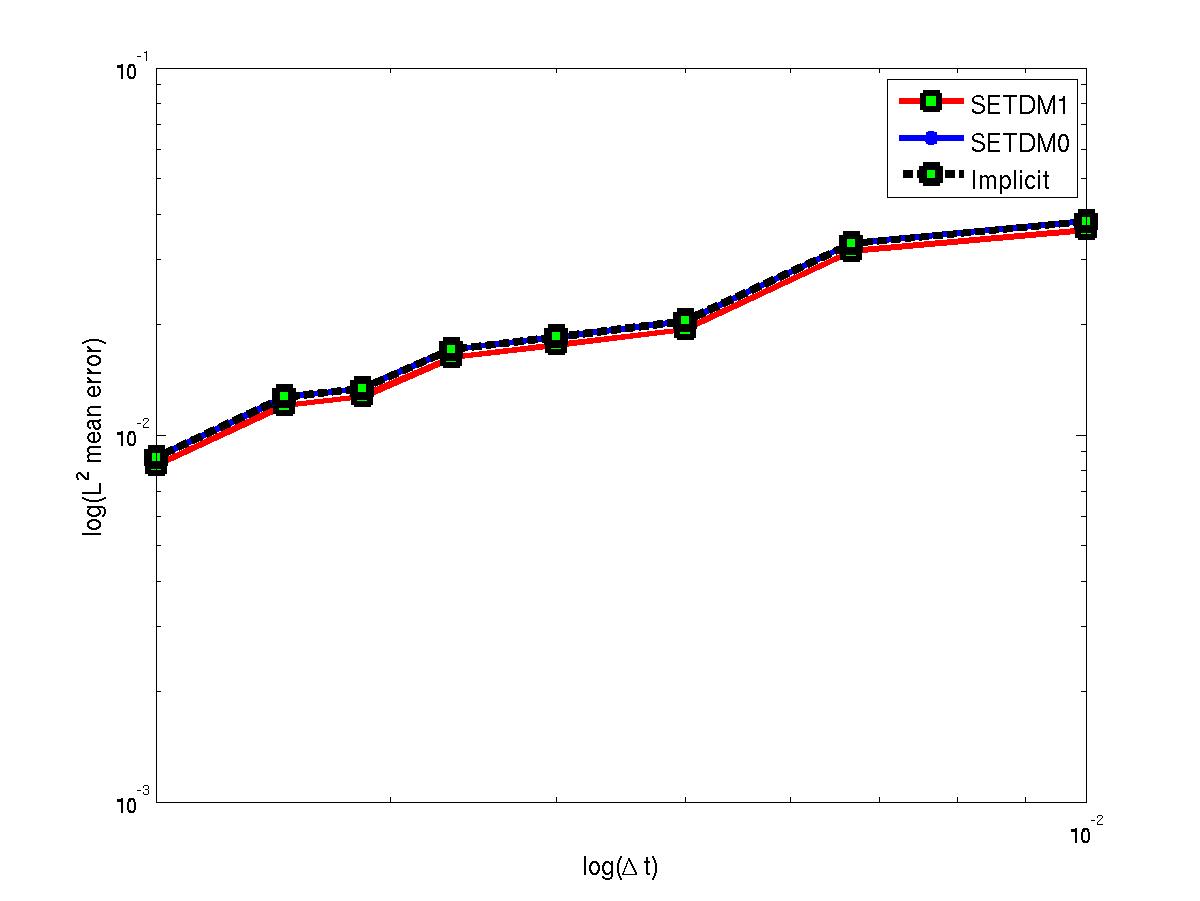}}
  \subfigure[]{
    \label{FIG0022b}
    \includegraphics[width=0.48\textwidth]{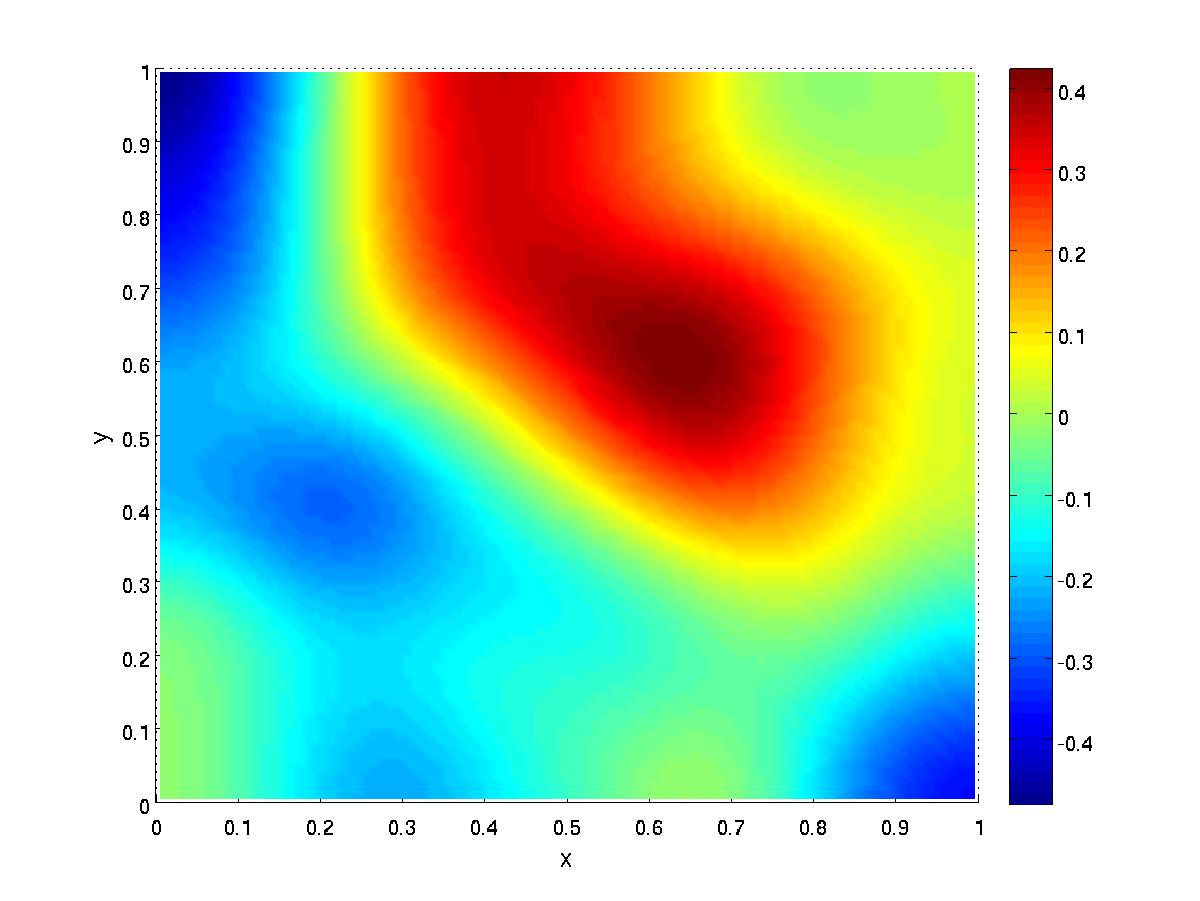}}
  \caption{ (a) Convergence of the root mean square $L^{2}$ norm at $T=1$ as
    a function of $\Dt$ with 10 realizations with $X_{0}=0,\;\Gamma=1,\; D=1$. The noise is white in time and
    with exponential correlation in space with lengths $b_1=b_2=0.2$.
    The temporal order of convergence in time is $0.9$ for all schemes. In (b) we plot a sample true solution.}    
  \label {FIG0022}
\end{figure}

We consider the covariance operator $Q$ with the following  covariance function (kernel) which has strong exponential decay 
\begin{eqnarray*}
 C_{r}((x_{1},y_{1});(x_{2},y_{2}))=\dfrac{\Gamma}{4 b_{1}b_{2}} \exp \left(-\dfrac{\pi}{4}\left[\dfrac{\left( x_{2}-x_{1}\right)^{2} }{b_{1}^{2}}+ \dfrac{\left( y_{2}-y_{1}\right)^{2} }{b_{2}^{2}}\right] \right),
\end{eqnarray*}
where $b_{1},b_{2}$ are  spatial correlation lengths in $x-$ axis and y- axis respectively and  $\Gamma>0$.\\
It is well known that the eigenfunctions  $ \{e_{i}^{(1)}e_{j}^{(2)}\}_{i,j\geq 0} $ of the operator $A=D \varDelta  $\;is given by
\begin{eqnarray}
\label{eigf}
\left\lbrace \begin{array}{l}
e_{0}^{(l)}=\sqrt{\dfrac{1}{L_{l}}},\;\;\;\lambda_{0}^{(l)}=0,\;\;\; e_{i}^{(l)}=\sqrt{\dfrac{2}{L_{l}}}\cos(\lambda_{i}^{(l)}x),\;\;\;\lambda_{i}^{(l)}=\dfrac{i \,\pi }{L_{l}}\\
\newline\\
l \in \left\lbrace 1, 2 \right\rbrace \;\;\; i=1, 2, 3, \cdots
 \end{array}\right.
\end{eqnarray}
with the corresponding eigenvalues $ \{\lambda_{i,j}\}_{i,j\geq 0} $ given by 
\begin{eqnarray*}
\lambda_{i,j}= (\lambda_{i}^{(1)})^{2}+ (\lambda_{j}^{(2)})^{2}.
\end{eqnarray*}
The corresponding values of $\left\lbrace q_{i,j}\right\rbrace_{i+j >
  0}$ in the representation (\ref{eq:W}) are given by  
\begin{eqnarray*}
 q_{i,j}= \Gamma \exp\left[ -\dfrac{1}{2 \pi}\left((\lambda_{i}^{(1)}b_{1})^{2}+(\lambda_{j}^{(2)}b_{2})^{2}\right) \right],
\end{eqnarray*}
see \cite{ATthesis} for details and
\cite{shardlow05,GrcaOjlvoSncho}. 
We compute the exponential functions $\varphi_{i},\;i=0,1$ with the real fast L\'{e}ja point technique and the absolute tolerance $10^{-6}$.
In our simulation  we take $L_{1}=L_{2}=1$ and the finite element triangulation is contructed with the rectangular 
grid with size $\Delta x= \Delta y = 1/150$. \figref{FIG0022a} shows the time convergence of SETDM1,SETDM0 and semi-implict schemes. 
The three methods have the same order of accuracy.
The temporal order of convergence that we observe is $0.9$ for all the schemes. This is higher than the predicted theoretical order of comvergence $0.5$ 
 in \thmref{th1}. The noise is regular and this order agrees with that in \cite{LR}.

%
As a more challenging example we consider the stochastic advection
diffusion reaction SPDE with multiplicative noise
\begin{eqnarray}
\label{advection}
dX&=&\left(\nabla \cdot  \mathbf{D}\nabla X -  \nabla \cdot (\mathbf{q}
  X)-\frac{X}{\vert X \vert+1}\right)dt+ X dW,\qquad  X_{0}=0 \qquad \Omega=[0,1]\times [0,1]\\
 \mathbf{D}&=&\left( \begin{array}{cc}
             10^{-2}&0\\
             0& 10^{-3}
           \end{array}\right)
\end{eqnarray}
with mixed Neumann-Dirichlet boundary conditions. The Dirichlet boundary condition is $X=1$ at $x=0$ and we use the homogeneous
Neumann boundary conditions elsewhere.
 According to \thmref{regsoluion}, we need to take the initial data 

\begin{figure}[h!]
   \begin{center}
  \subfigure[]{
    \label{FIG022a}
    \includegraphics[width=0.32\textwidth]{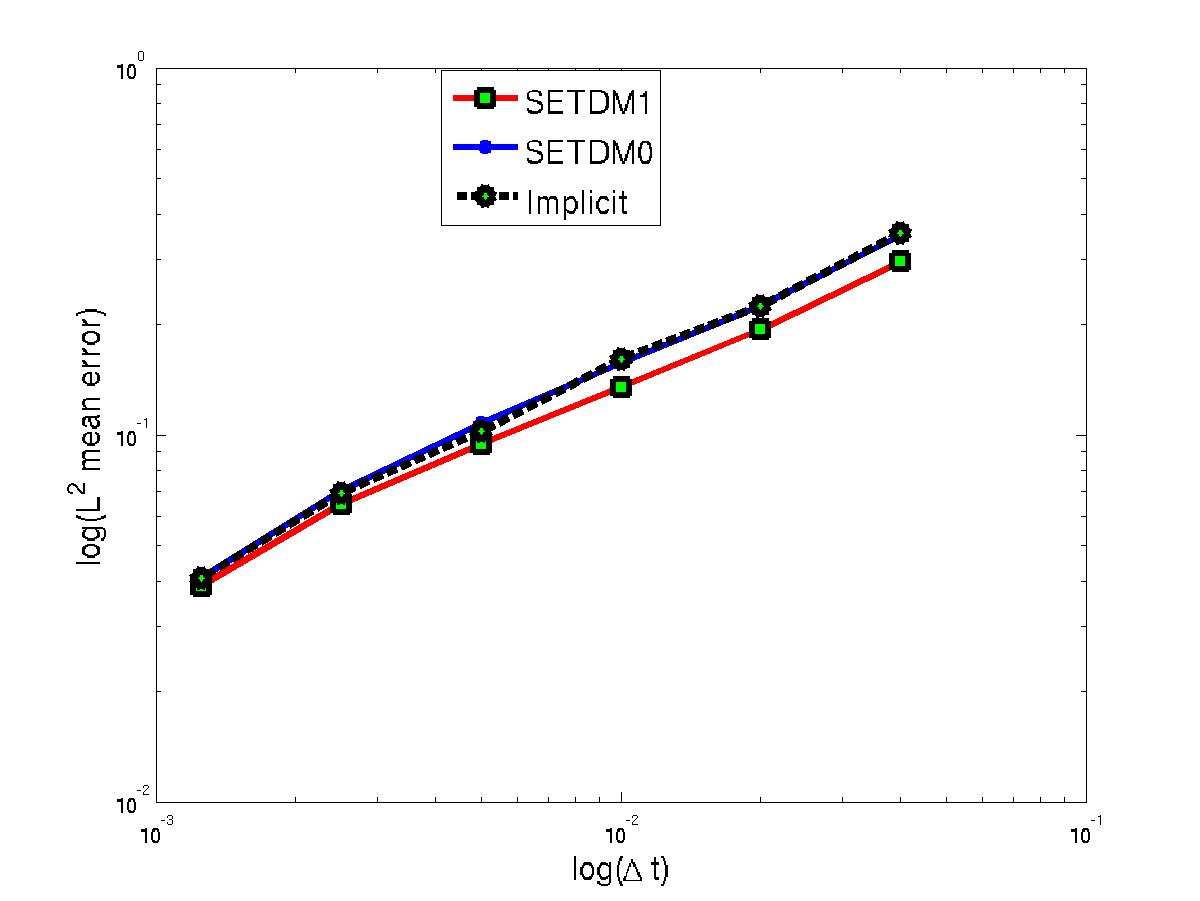}}
  \subfigure[]{
    \label{FIG022b}
    \includegraphics[width=0.32\textwidth]{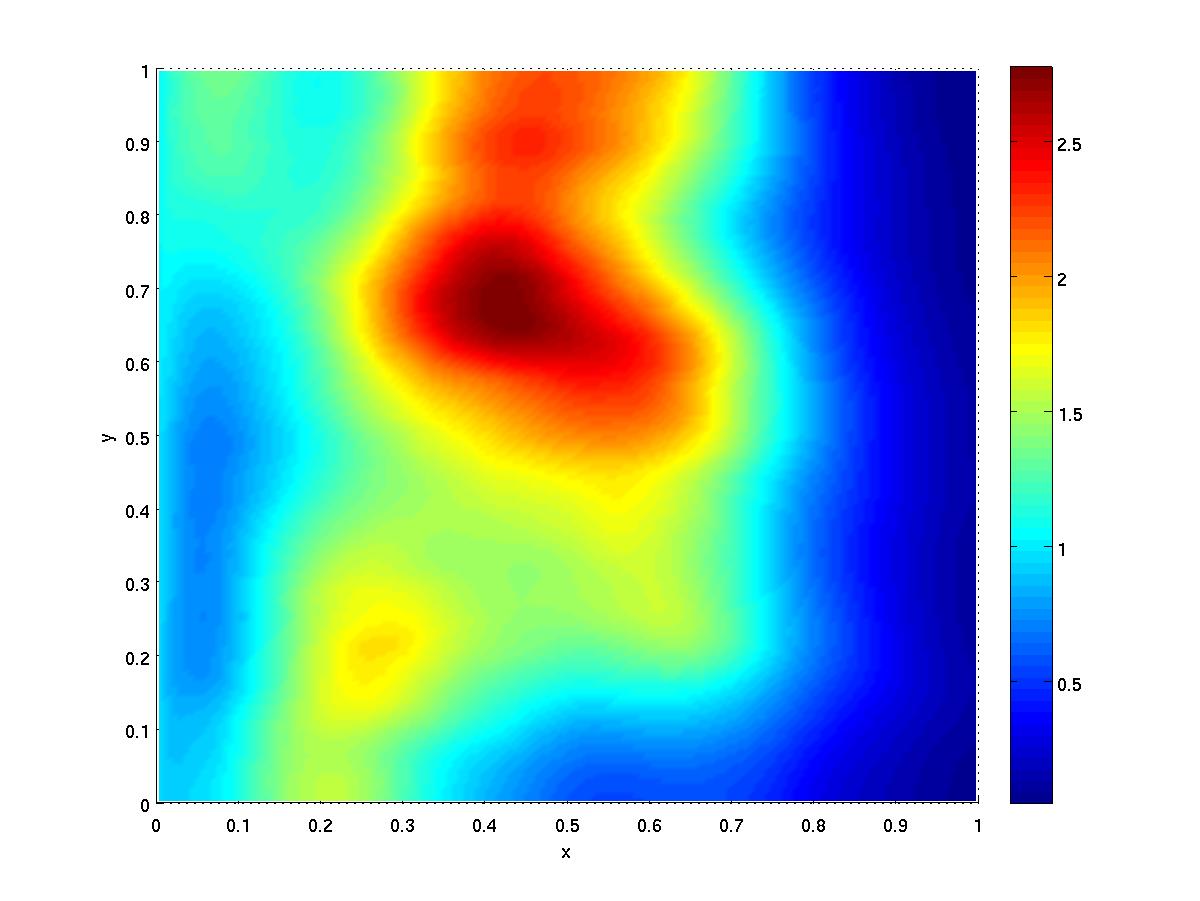}}
   \subfigure[]{
    \label{FIG022c}
    \includegraphics[width=0.32\textwidth]{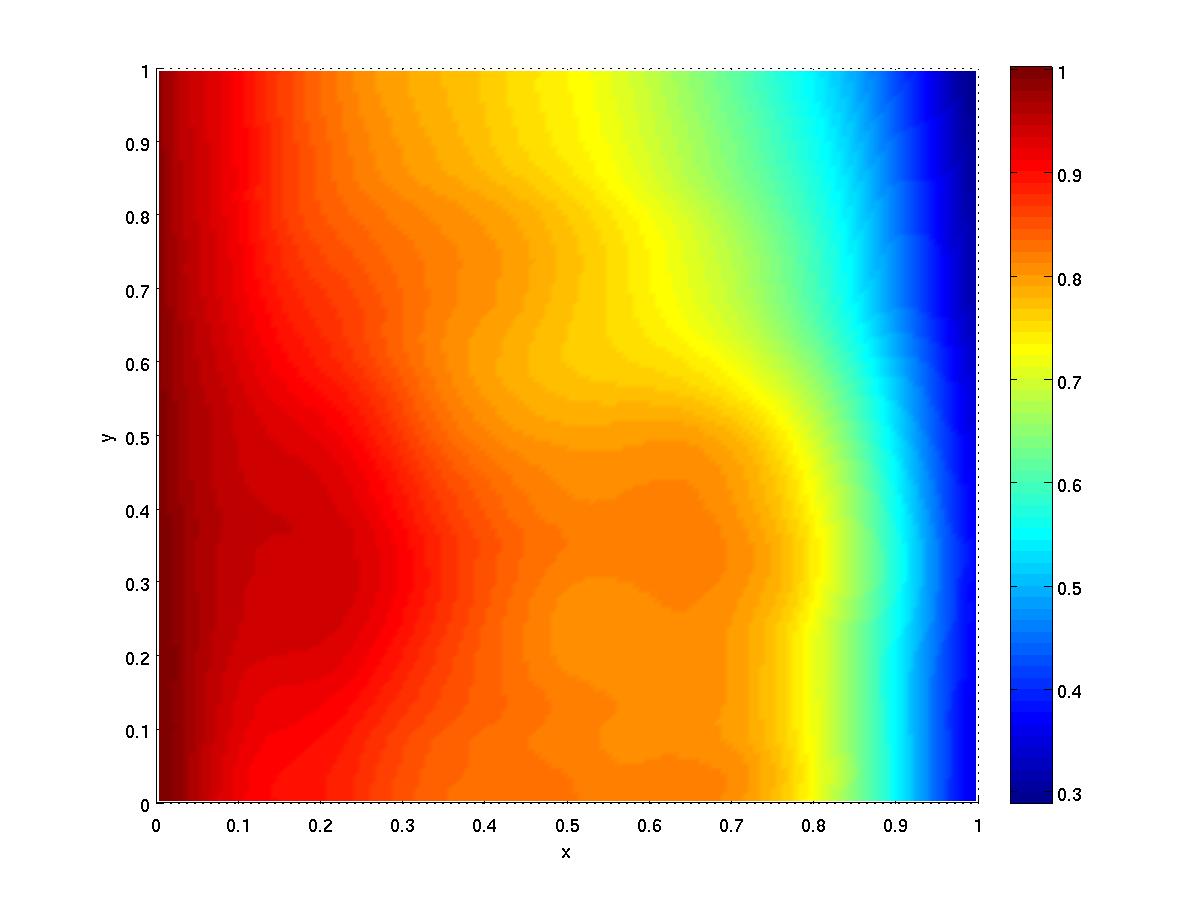}} 
  \caption{ (a) Convergence of the root mean square $L^{2}$ norm at $T=1$ as
    a function of $\Dt$ with 200 realizations with $\Delta x= \Delta y
    = 1/150$, $X_{0}=0$ for the homogeneous medium. The
    noise is white in time with \eqref{noise2} and $r=2.01$.  
    The temporal orders of convergence in time are $0.54$ for SETDM1 and
    $0.58$ for SETDM0 semi-implicit schemes. In
    (b) we plot a sample of the 'true solution' for $r=2.01 $ with 
    $\Dt=1/1600$ while (c) shows the mean of the “true solution“ for 200
realizations.}
  \label {FIG022}
\end{center}
\end{figure}

\begin{figure}[h!]
  \begin{center}
  \subfigure[]{
    \label{FIG027a}
    \includegraphics[width=0.39\textwidth]{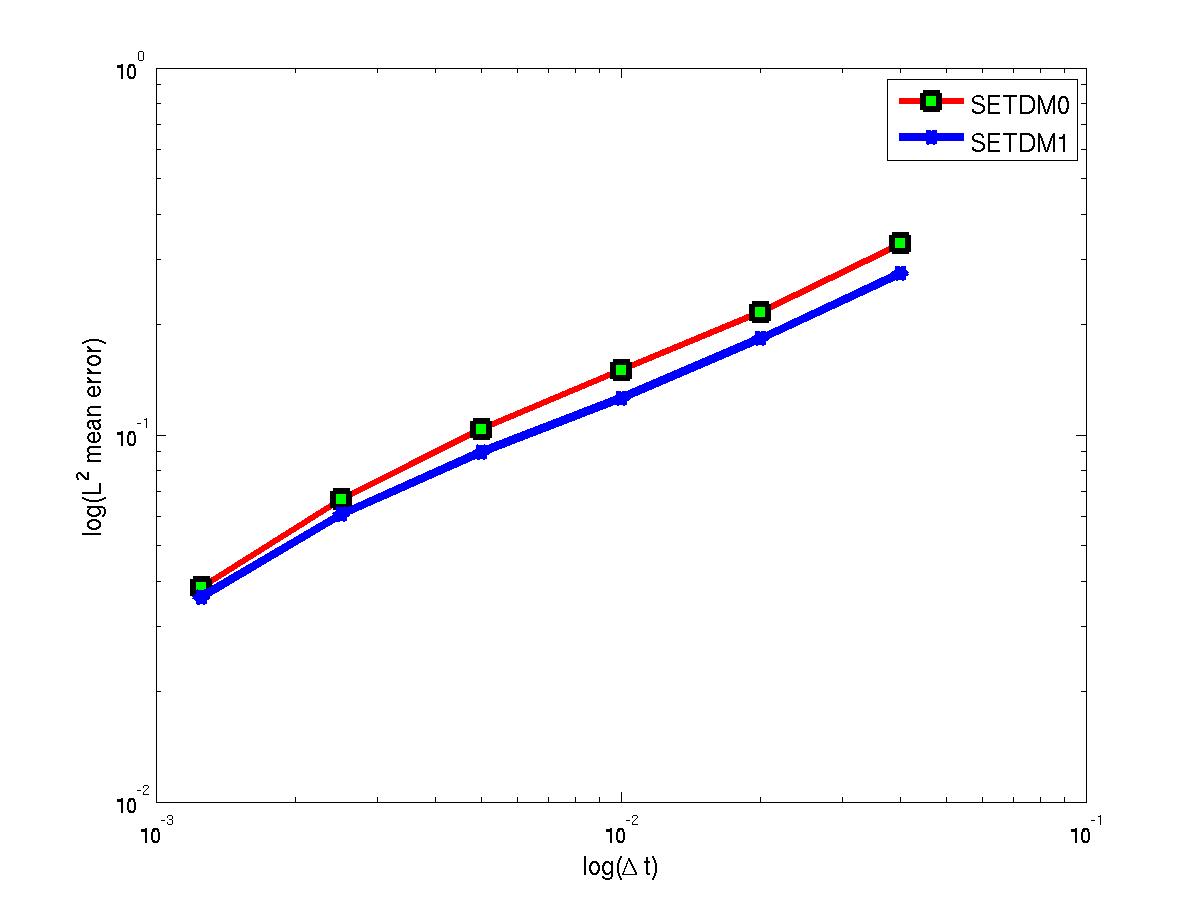}}
  \subfigure[]{
    \label{FIG027b}
    \includegraphics[width=0.39\textwidth]{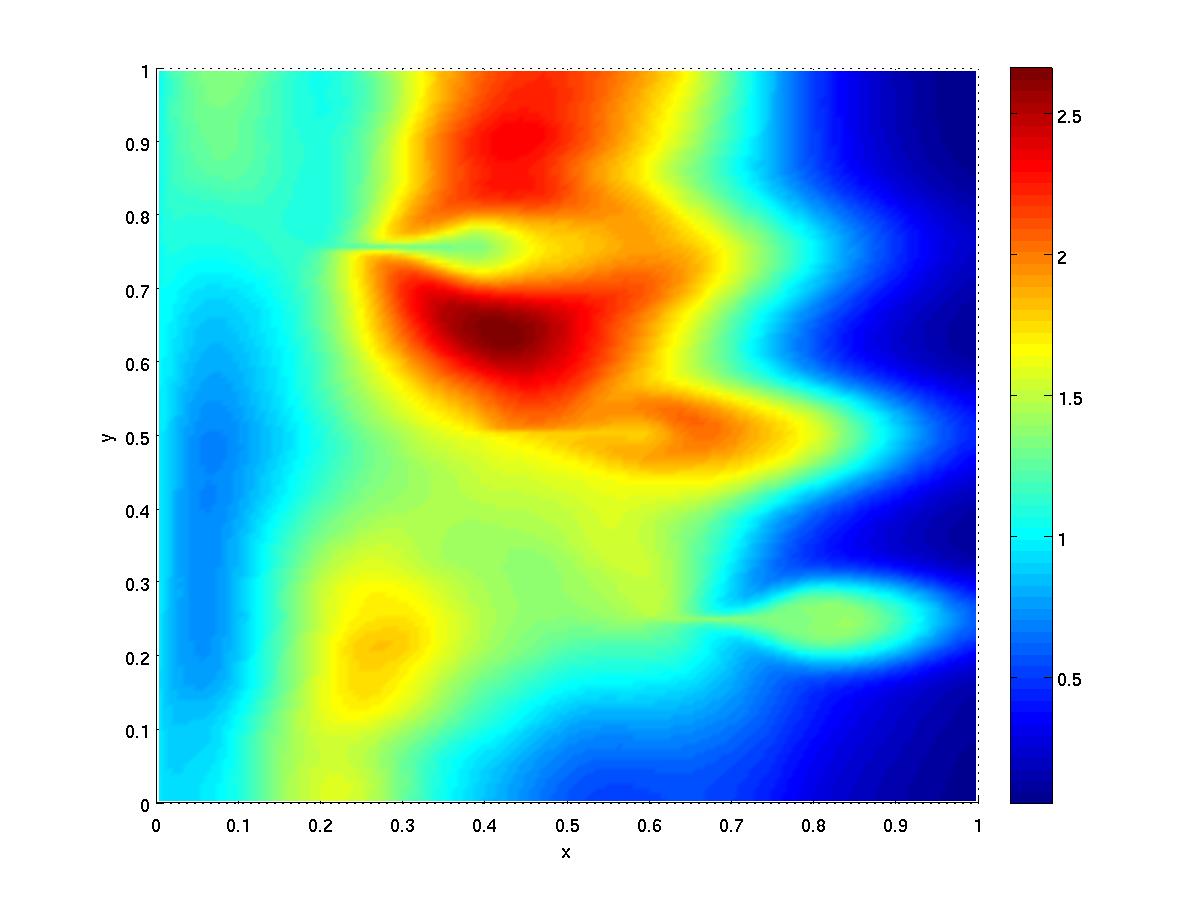}}
\subfigure[]{
    \label{FIG027c}
    \includegraphics[width=0.39\textwidth]{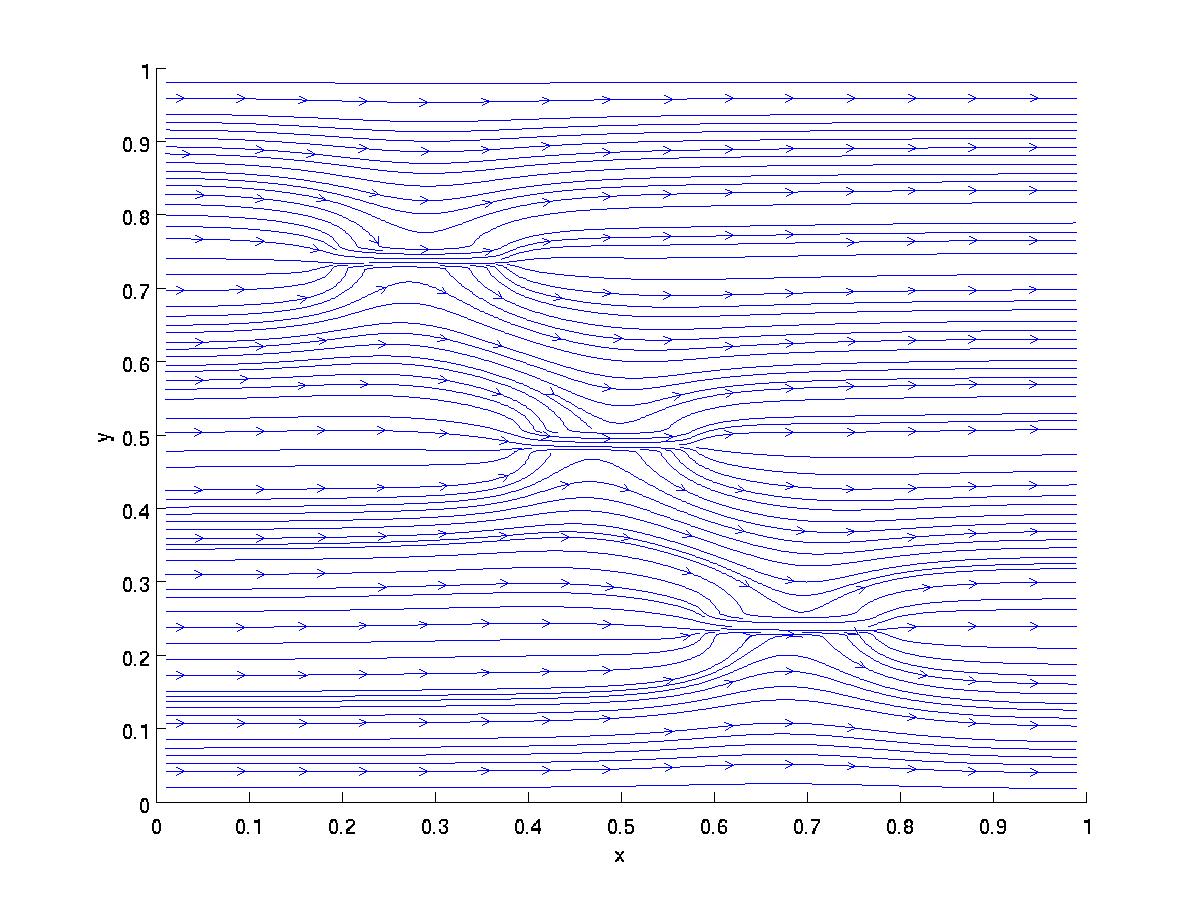}}
\subfigure[]{
    \label{FIG027d}
    \includegraphics[width=0.39\textwidth]{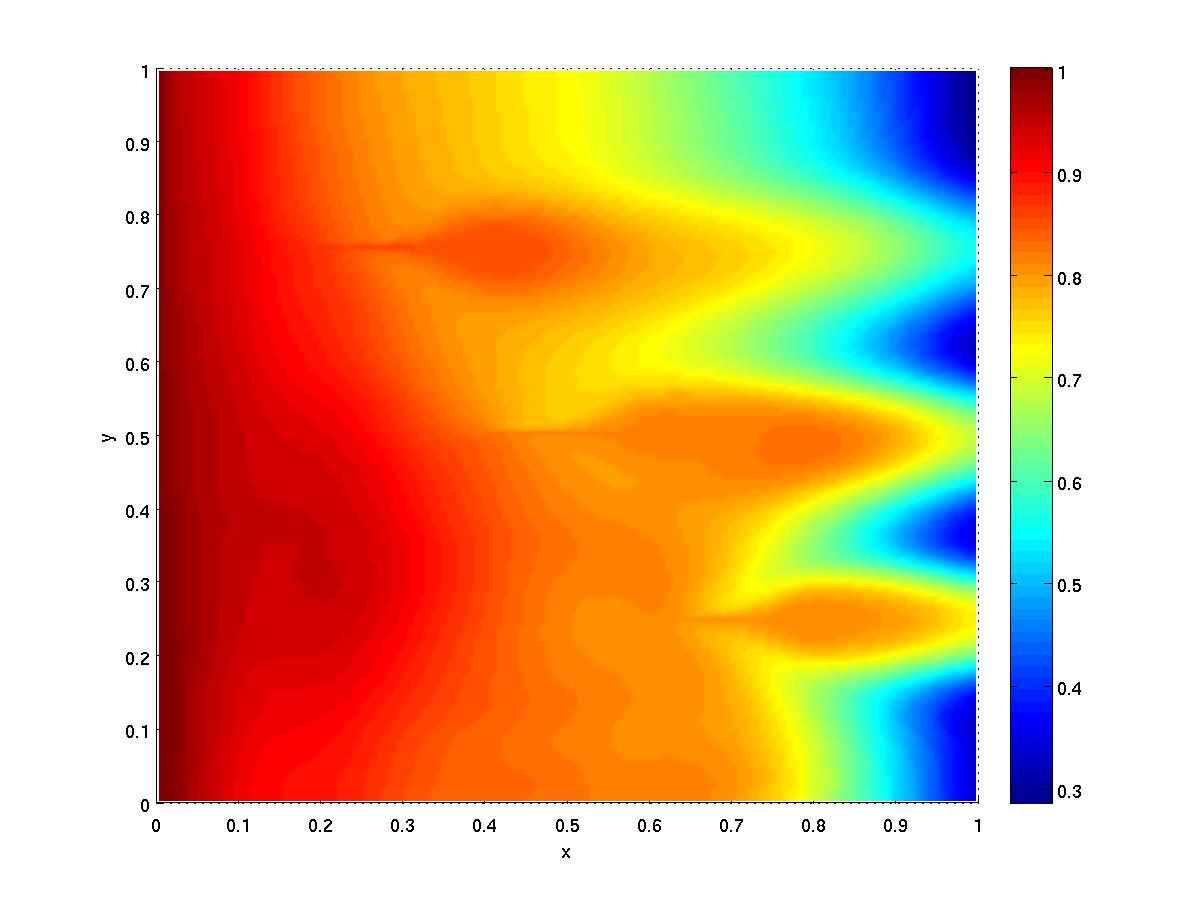}}       
  \caption{ (a) Convergence of the root mean square $L^{2}$ norm at $T=1$ as
    a function of $\Dt$ with 200 realizations with $\Delta x= \Delta y
    = 1/150$, $X_{0}=0$ for the heterogeneous medium.  The noise is white in time with \eqref{noise2} and $r=2.01$.
    The temporal orders of convergence in time are $0.54$ for  SETDM1 scheme and $0.58 $ for SETDM0. In (b) we plot a sample 'true solution' for $r=2.01$ with
    $\Dt=1/1600$. In (c) we plot the streamline of the velocity field while (d) shows the mean of the ``true solution`` for 200 realizations.}    
  \label {FIG072}
 \end{center}
\end{figure}

$X_{0} \in L_{2}(\mathbb{D},\mathcal{D}((-A)^{\beta})),\; \beta >0$ to have a regular solution such that
$XdW$ make sense. For our simulation we take $X_{0}=0$.
For a  homogeneous medium, we use the constant velocity $\mathbf{q}=(1,0)$. 
In terms of equation \eqref{sadr} the nonlinear terms $f$ and $b$ are given by
\begin{eqnarray}
 f(x,u)=-\frac{u}{(\vert u \vert +1)}, \qquad  b(x,u) =u, \quad u \in
 \mathbb{R},\; x \in \Omega,
\end{eqnarray}
and the corresponding Nemytskii  operators $F$  and $B$ are obtained from \eqref{nemform}
and  clearly satisfy \assref{assumption1} (if the domain of $f$ is restricted to $\mathbb{R}^{+}$) 
and \assref{assumption2} (see \cite[Section 4]{MJentzen1}) respectively, where  \eqref{eigf} is used in the noise representation (\ref{eq:W}).
The linear operator is  given by
\begin{eqnarray}
 A=\nabla \cdot  \mathbf{D}\nabla (.) -  \nabla \cdot \mathbf{q}(.).
\end{eqnarray}
%
For a heterogeneous medium we used three parallel high 
 permeability streaks. This could represent for example a highly idealized fracture pattern.
  We obtain the Darcy velocity field $\mathbf{q}$  by  solving  the system
\begin{eqnarray}
\label{darcy}
 \left\lbrace \begin{array}{l}
 \nabla \cdot\mathbf{q} =0\\
 \mathbf{q}=-\dfrac{k(\mathbf{x})}{\mu} \nabla p,
\end{array}\right.
\end{eqnarray}
with  Dirichlet boundary conditions 
$\Gamma_{D}^{1}=\left\lbrace 0 ,1 \right\rbrace \times \left[
  0,1\right] \text{and Neumann boundary} 
 \quad \Gamma_{N}^{1}=\left( 0,1\right)\times\left\lbrace 0
   ,1\right\rbrace $ such that 
\begin{eqnarray*}
 p&=&\left\lbrace \begin{array}{l}
1 \quad \text{in}\quad \left\lbrace 0 \right\rbrace \times\left[ 0,1\right]\\
0 \quad \text{in}\quad \left\lbrace 1 \right\rbrace \times\left[ 0,1\right]
 \end{array}\right. 
\end{eqnarray*}
and 
$$
- k \,\nabla p (\mathbf{x},t)\,\cdot \mathbf{n} =0\quad \text{in}\quad \Gamma_{N}^{1},
$$
 where  $p$ is the pressure, $\mu$ is dynamical viscosity and $k$ the
 permeability of the porous medium.  
 We have assumed that rock and fluids are incompressible and sources or
 sinks are absent, thus the equation 
\begin{eqnarray}
  \label{pr}
  \nabla \cdot\mathbf{q}=\nabla \cdot 
  \left[\dfrac{k(\mathbf{x})}{\mu} \nabla p\right]=0 
\end{eqnarray}
comes from mass conservation.
 As in  \cite{LT,KLNS}, we take the following values for
$\left\lbrace q_{i,j}\right\rbrace_{i+j > 0}$ in the representation 
(\ref{eq:W}) 
\begin{eqnarray}
  \label{noise2}
  q_{i,j}= 1/\left( i+j\right)^{r},\qquad r>0.
\end{eqnarray}
Note that to have a trace class noise we need $r>2$. In our simulation we use  $r=2.01$. 
To deal with high  P\'{e}clet flows we discretize in space using
finite volumes. We can write the semi-discrete finite volume discretization 
of \eqref{advection} as 
\begin{eqnarray}
 dX^{h}=(A_{h}X^{h}+P_{h}F(X^{h})) + P_{h} B(X^{h})dW,
\end{eqnarray}
(see \cite{FV,ATthesis}).
\figref{FIG022a} shows the convergence of SETDM0, SETDM1 and
semi-implicit schemes for the homogeneous porous medium. 
The scheme SETDM1 seems to be more accurate
 for large time steps but for large time steps it has the same order of accuracy as the semi-implicit  and SETDM0 schemes. The temporal order is
 $0.54$ for  SETDM1 scheme, $0.58$ for SETDM0 scheme 
and the semi-implicit scheme. We used 200 realizations and the convergence order is close to the $0.5$, the predicted order of convergence in \thmref{th1}.
 A sample the 'true solution' is shown in  \figref{FIG022b} with $\Dt=1/1600$ while the mean of the ''true solution'' for 200
realizations is shown in \figref{FIG022c}. 
\newpage
\figref{FIG027a} shows the convergence of SETDM0 and SETDM1 schemes
for the heterogeneous porous medium. It also shows that SETDM1 
is more accurate than SETDM0 scheme for high time step size. The observed temporal order is
 $0.54$ for  SETDM1 scheme and $0.58$ for SETDM0 scheme. \figref{FIG027c}
shows the streamline of the velocity field. 
A sample the ``true solution`` is shown in  \figref{FIG027b} with
$\Dt=1/1600$ while the mean of the ''true solution'' for 200
realizations is shown in  \figref{FIG027d}. 

To conclude we have proved the errors estimates for the exponential
based integrators  and observed the predicted rate of convergence in
the simulations.


\end{document}